\documentclass{amsart}

\usepackage{color}
\usepackage{epsfig}
\usepackage{amsthm}
\usepackage{amssymb}
\usepackage{amsmath}
\usepackage{amscd}
\usepackage{yhmath}

\usepackage{framed}

\newtheorem {Theorem}  {Theorem}

\numberwithin{Theorem}{section}

\newtheorem {Lemma}[Theorem]  {Lemma}
\newtheorem {Proposition}[Theorem]{Proposition}
\theoremstyle{definition}
\newtheorem{Definition}[Theorem]{Definition}
\theoremstyle{remark}
\newtheorem{Remark}[Theorem]{Remark}
\newtheorem{Example}[Theorem]{Example}

\newtheorem {Corollary}[Theorem]{Corollary}

\expandafter\chardef\csname pre amssym.def
at\endcsname=\the\catcode`\@ \catcode`\@=11
\def\undefine#1{\let#1\undefined}
\def\newsymbol#1#2#3#4#5{\let\next@\relax
 \ifnum#2=\@ne\let\next@\msafam@\else
 \ifnum#2=\tw@\let\next@\msbfam@\fi\fi
 \mathchardef#1="#3\next@#4#5}
\def\mathhexbox@#1#2#3{\relax
 \ifmmode\mathpalette{}{\m@th\mathchar"#1#2#3}%
 \else\leavevmode\hbox{$\m@th\mathchar"#1#2#3$}\fi}
\def\hexnumber@#1{\ifcase#1 0\or 1\or 2\or 3\or 4\or 5\or 6\or 7\or 8\or
 9\or A\or B\or C\or D\or E\or F\fi}

\font\teneufm=eufm10 \font\seveneufm=eufm7 \font\fiveeufm=eufm5
\newfam\eufmfam
\textfont\eufmfam=\teneufm \scriptfont\eufmfam=\seveneufm
\scriptscriptfont\eufmfam=\fiveeufm

\catcode`\@=\csname pre amssym.def at\endcsname

\newcounter{remark}
\def  \eps  {\epsilon}

\def\3n{\negthinspace \negthinspace \negthinspace }
\def\2n{\negthinspace \negthinspace }
\def\1n{\negthinspace }

\newcommand{\bg}{\begin{equation}}
\newcommand{\ed}{\end{equation}}
\newcommand{\bga}{\begin{eqnarray}}
\newcommand{\eda}{\end{eqnarray}}
\newcommand{\pf}{\textbf{Proof \ }}

\def\cbdu{\hfill{$\Box$}}

\newcommand{\ueps}{u^{\eps}}
\newcommand{\e}{\epsilon}

\renewcommand{\a}{\alpha}

\renewcommand{\b}{\beta}

\newcommand{\R}{\mathbf{R}}

\newcommand{\Mm}{{\mathcal M}}

\newcommand{\eledos}{L^2}

\def  \R   {{\mathbb R}}

\def  \12  {{\frac{1}{2}}}

\def\bd{\begin{Definition}}
\def\ede{\end{Definition}}
\def\be{\begin{equation}}
\def\bel{\begin{equation}\label}
\def\ee{\end{equation}}
\def\bt{\begin{Theorem}}
\def\et{\end{Theorem}}
\def\bex{\begin{Example}}
\def\eex{\end{Example}}
\def\bc{\begin{Corollary}}
\def\ec{\end{Corollary}}
\def\bl{\begin{Lemma}}
\def\el{\end{Lemma}}
\def\bp{\begin{Proposition}}
\def\ep{\end{Proposition}}

\def\br{\begin{Remark}}
\def\er{\end{Remark}}
\def\ba{\begin{array}}
\def\ea{\end{array}}

\def\bea{\begin{eqnarray}}
\def\eea{\end{eqnarray}}

\begin{document}

\title[Decay characterization for dissipative equations]{Decay characterization of solutions to dissipative equations}

\author{C\'esar J. Niche}
\address[C.J. Niche]{Departamento de Matem\'atica Aplicada, Instituto de Matem\'atica, Universidade Federal do Rio de Janeiro, CEP 21941-909, Rio de Janeiro - RJ, Brasil}
\email{cniche@im.ufrj.br}

\author{Mar\'{\i}a E. Schonbek}
\address[M.E. Schonbek]{Department of Mathematics, UC Santa Cruz, Santa Cruz, CA 95064, USA}
\email{schonbek@ucsc.edu}

\thanks{C.J. Niche acknowledges financial support from PRONEX  E-26/110.560/2010-APQ1, FAPERJ-CNPq and Ci\^encia sem Fronteiras - PVE 011/12. \\
M. E. Schonbek was partially supported by NSF Grant DMS-0900909.}

\begin{abstract}
We address the study of decay rates of solutions to dissipative equations. The characterization of these rates is given for a wide class of  linear systems by the {\em decay character}, which is a number  associated to the  initial datum that  describes the behavior of the datum near the origin in frequency space. We then use the decay character and the Fourier Splitting method to obtain upper and lower bounds for decay of solutions to appropriate dissipative nonlinear equations, both in the incompressible and compressible case. 
\end{abstract}

\date{\today}

\maketitle

\section{Introduction}

Solutions to many linear and nonlinear dissipative evolution equations obey inequalities of the form

\be
\label{eqn:ineq-fourier-splitting}
\frac{1}{2} \frac{d}{dt} \Vert f(t) \Vert _{L^2} ^2 \leq  - C  \int _{\R^n} |\xi|^{2 \a} |\widehat{f} (\xi, t)|^2 \, d \xi, \qquad 0 < \a \leq 1, \, C > 0
\ee
as  for example is the case for the Navier-Stokes equations, where $\a = 1$ and $f$ is the velocity field, or for  the $2D$ quasi-geostrophic equation, where $0 < \a \leq 1$ and $f$ is the scalar potential temperature. Inequality (\ref{eqn:ineq-fourier-splitting}) indicates that the $L^2$ norm of solutions decays in time, hence  it is natural to ask what is the decay rate for these solutions.

In order to provide  context to our work, we recall  some well-known results concerning the decay of solutions to the initial value problem for the $3D$ Navier-Stokes equations

\begin{eqnarray}
\label{eqn:navier-stokes}
\partial _t u + (u \cdot \nabla) u & = & \Delta u - \nabla p, \nonumber \\ div \, u & = & 0, \nonumber \\ u_0 (x) & = & u(x,0). 
\end{eqnarray}
 When the initial datum $u_0$ is in $L^2$, Masuda \cite{masuda1984} showed that the $L^2$ norm of weak solutions goes to zero as time goes to infinity.  If $u_0 \in L^1 \cap L^2$, M.E. Schonbek \cite{MR775190}, \cite{MR837929} proved that the decay  has a uniform rate

\begin{displaymath}
\Vert u(t) \Vert _{L^2} \leq C (1 + t) ^{- \frac{3}{4}}, \qquad t > 0.
\end{displaymath}
We note that if the initial datum is small, the above result had been obtained by Kato \cite{Kato1984}.
 M.E. Schonbek \cite{MR837929} also proved that if $u_0$ is just in $L^2$, there are solutions, whose initial data is a rescaled Gaussian, that do not have a uniform algebraic decay rate. Later, Wiegner \cite{MR881519}, for $2 \leq n \leq 4,$ proved more precise uniform decay rates by carefully analyzing the relation between the linear and nonlinear part of the solution: if $\Vert e^{t \Delta} u_0 \Vert _{L^2} \leq C (1 + t)^{- \mu}$, for $\mu \geq 0$, then

\begin{displaymath}
\Vert u(t) \Vert _{L^2} \leq C (1 + t) ^{- \min \{ \mu, \frac{5}{4} \}}, \qquad t > 0.
\end{displaymath}
Note that when $u_0 \in L^1 \cap L^2$, we have that $\widehat{u}_0 \in L^{\infty}$, so $\Vert e^{t \Delta} u_0 \Vert _{L^2} \leq C (1 + t)^{- \frac{3}{4}}$, which leads to the previous estimate. The faster decay rate is obtained, for example, for solutions with zero average $u_0$.

Recently, Bjorland and M.E. Schonbek \cite{MR2493562} introduced a new idea in this area, by associating to every  $u_0 \in L^2$ a {\em decay character} $r^{\ast} = r^{\ast} (u_0)$, which describes the rate with which the solution to (\ref{eqn:navier-stokes}) with such initial datum decays. A key point in the proof of their results is an estimate for the decay rate of solutions to the heat equation in terms of $r^{\ast} (u_0)$.

The main goal of this article is to refine and extend Bjorland and M.E.Schonbek's work. First, we define the decay character of $\Lambda^s u_0$, for $u_0 \in H^s (\R^n)$, $(- \Delta) ^{\frac{1}{2}} = \Lambda$ and $s \geq 0$ and  establish its relation with the decay character of $u_0$. Then, we give sharp upper and lower bounds for decay rates in Sobolev space for a wide class of linear equations in terms of the decay characters of $u_0$ and $\Lambda^s u_0$. Finally, we use these estimates to obtain results for decay of solutions to some nonlinear dissipative equations.

We next describe in detail the definitions introduced and the results obtained in this article. It has been frequently remarked that {\em the long time behavior in time of solutions is determined by the small frequencies of the initial data/solution}. In order to formalize this, we want to know what is the order of $\widehat{\Lambda ^s u}_0 (\xi)$ at the origin, by comparing it to $f(\xi) = |\xi|^{r}$. Thus, in Section \ref{section-decay-character} we define  the {\em decay indicator} of $\Lambda^s u_0$, for  every $u_0 \in L^2 (\R^n)$, as

\begin{displaymath}
P_r^s(u_0) = \lim _{\rho \to 0} \rho ^{-2r-n} \int _{B(\rho)} \large| \xi|^{2s}|\widehat{u}_0 (\xi) \large|^2 \, d \xi
\end{displaymath} 
where $s \geq 0$, $\Lambda = (- \Delta) ^{\frac{1}{2}}$ and $B(\rho)$ is the ball at the origin with radius $\rho$.  By taking $r = q + s$, we see that the decay indicator compares $|\widehat{\Lambda^s u_0} (\xi)|^2$ to $f(\xi) = |\xi|^{2(q + s)}$ at $\xi = 0$. When $s = 0$, we recover Definition 5.5 from Bjorland and M.E. Schonbek \cite{MR2493562} and we use the notation $r^{\ast} _0 = r^{\ast}$. If for some finite $r$ we have that $\widehat{\Lambda^s u}_0 (\xi)$ and $|\xi|^{2(q + s)}$ are equivalent, then $P_r^s(u_0) > 0$ and we say the {\em decay character} of $\Lambda^s u_0$, is $r_s^{\ast} = r_s^{\ast}( u_0) = r$. If $\widehat{\Lambda^s u}_0 (\xi)$ goes to zero at the origin faster (respectively slower) than any $|\xi|^{2(q + s)}$, we say that the {\em decay character} is $r_s^{\ast} = \infty$ (respectively $r_s^{\ast} = - \frac{n}{2}+s$).

Inequality (\ref{eqn:ineq-fourier-splitting}) is key for obtaining decay for both linear and nonlinear dissipative equations. In Section \ref{section-linear-part} we introduce a family of pseudodifferential operators $\mathcal{L}$ that are negative definite and diagonalizable. This family contains, amongst others, the usual Laplacian and fractional Laplacian, as well as a compressible approximation to the Stokes operator which is the linear part of one of the  nonlinear equations we analize later. For all these operators, the linear equation $v_t = \mathcal{L} v$ leads to (\ref{eqn:ineq-fourier-splitting}). 

In Sections \ref{section-decay-l2-linear-part} and \ref{section-decay-hs-solution} we establish the main results concerning  decay of linear equations in terms of the decay character. In Theorem \ref{characterization-decay-l2} in Section \ref{section-decay-l2-linear-part}, the decay of the $L^2$ norm of solutions is characterized for data  $ v_0 \in L^2$. More specifically,  depending on the decay character $r^{\ast} = r^{\ast} (v_0)$, decay is faster or slower than any algebraic rate or is bounded from above and below by algebraic rates with the same exponent. Section \ref{section-decay-hs-solution}  deals with the characterization of decay of the homogeneous Sobolev norm $\dot{H}^s$ of solutions, when $v_0\in H^s$.  Theorem \ref{decay-character-hs}  establishes a relation between the decay characters  
$r^{\ast}$ and $r_s ^{\ast} = r_s ^{\ast} (v_0)$,   proving that $r_s ^{\ast} = s + r ^{\ast}$. As $r_s ^{\ast}$ is the decay character of the $L^2$ function $\Lambda ^s v_0$, the characterization of decay for the  $\dot{H}^s$ norm in Theorem \ref{characterization-decay-hs}, analogous to the one obtained in Theorem \ref{characterization-decay-l2} for $v_0$ in $L^2 (\R^n)$, follows from Theorems  \ref{characterization-decay-l2} and \ref{decay-character-hs}.

In Section \ref{applications} we  study  decay rates of solutions to nonlinear equations, using the  decay results of their linear part. We first consider  the dissipative quasi-geostrophic equation  in $\R^2$

\begin{displaymath}
\theta _t + u \cdot \nabla \theta + \kappa (- \Delta) ^{\a}u = 0, \qquad 0 < \a \leq 1, \kappa > 0
\end{displaymath}
where $u = R^{\bot} \theta = \left( -R_2 \theta, R_1 \theta \right)$, $R_i$ is the Riesz transform in the $x_i$ variable and $\theta = \theta (x,t)$ is a scalar. This active scalar equation has been extensively studied in the last 20 years because for $\a = \frac{1}{2}$, provides a dimensionally correct $2D$ model to the $3D$ Navier-Stokes equations. Moreover, this equation and its inviscid counterpart (i.e., $\kappa = 0$), describe important models of currents and atmosphere circulation. Using the Fourier Splitting method, developed by M.E. Schonbek \cite{MR571048}, \cite{MR775190}, \cite{MR837929} to study decay of solutions to conservation laws and the Navier-Stokes equations, in Theorems \ref{decay-l2-qg} and \ref{lower-bounds-qg} we prove upper and lower bounds for the decay rate of this equation, to wit, we obtain that for $\theta_0$ in $L^2 (\R^2)$ with $r^{\ast} = r^{\ast} (\theta_0)$, we have that for $r^{\ast} \leq 1 - \a$ and some $C_1, C_2 > 0$

\be
\label{eqn:sharp-qg}
C_1 (1 + t) ^{-\frac{1}{\a} (1 + r^{\ast})} \leq \Vert \theta (t) \Vert _{L^2} ^2 \leq C_2 (1 + t) ^{-\frac{1}{\a} (1 + r^{\ast})};
\ee
in the region determined by $r^{\ast} \geq 1 - \a$, $r^{\ast} \leq 1$, $r^{\ast} \leq 2(1 - \a)$ and for some $C_1, C_2 > 0$ we have

\be
\label{eqn:gap-qg}
C_1 (1 + t) ^{-\frac{1}{\a} (1 + r^{\ast})} \leq \Vert \theta (t) \Vert _{L^2} ^2 \leq C_2 (1 + t) ^{-\frac{1}{\a} (2 - \a)},
\ee
and in the region determined by $r^{\ast} > 1$ and $r^{\ast} \geq 2(1 - \a)$ we have that

\be
\label{eqn:no-lower-qg}
\Vert \theta (t) \Vert _{L^2} ^2 \leq C_2 (1 + t) ^{-\frac{1}{\a} (2 - \a)}.
\ee
Then, (\ref{eqn:sharp-qg}) provides a sharp characterization of decay in terms of the decay character. However, in (\ref{eqn:gap-qg}) we may have a gap in the decay rates and in (\ref{eqn:no-lower-qg}) we have no lower bound at all. Note that these problems arise when the linear part has a relatively fast decay rate. This phenomenon also arises in the Navier-Stokes equations, see Theorem A in Miyakawa and M.E. Schonbek \cite{MR1844282} and Theorem 6.5 in Bjorland and M.E. Schonbek \cite{MR2493562}. In Theorem \ref{decay-derivative-qg} we prove results concerning the decay of the $\dot{H}^s$ norm of solutions. 

In Section \ref{applications}  we also study the compressible approximation to Navier-Stokes equations

\begin{eqnarray*}
\partial _t \ueps  + (\ueps \cdot \nabla) \ueps  + \frac{1}{2} (div \, \ueps) \ueps & = & \Delta \ueps + \frac{1}{\eps} \, \nabla \cdot div \, \ueps  \nonumber \\ \ueps (x,0)  & = & \ueps _0 (x),
\end{eqnarray*}
introduced by Temam \cite{MR0237972}. This system is obtained from the Navier-Stokes equations (\ref{eqn:navier-stokes}) by relating the pressure $p$ to the velocity $u$ through $\eps p = - div \, u$ in order to eliminate the nonlocal relation between them. The nonlinear damping term $\frac{1}{2} (div \, \ueps) \ueps$ has to be added to have an energy inequality. This, together with the fact that the operator which defines the linear part of this system fits in the context of Section \ref{section-linear-part}, leads to (\ref{eqn:ineq-fourier-splitting}), which allows us to use the Fourier Splitting method. We then prove, in Theorems \ref{decay-rusin-l2} and \ref{lower-bound-rusin}, upper and lower bounds for the decay rate of solutions. More precisely, for $\ueps _0 \in L^2 (\R^3)$ with $r^{\ast} = r^{\ast} (u_0)$, we have that  for $-\frac{3}{2} <  r^* \leq 1$, there exist $C_1, C_2 > 0$ such that

\begin{displaymath}
C_1 (1 + t) ^{- \left( \frac{3}{2} + r^{\ast} \right)} \leq \Vert \ueps (t) \Vert _{L^2} ^2 \leq C_2 (1 + t) ^{- \left( \frac{3}{2} + r^{\ast} \right)},
\end{displaymath}
while for $r^{\ast} > 1$, we obtain

\begin{displaymath}
\Vert \ueps (t) \Vert _{L^2} ^2 \leq C (1 + t) ^{- \frac{5}{2}}.
\end{displaymath}
As before, when the linear part of the solution has relatively slow decay, we can sharply estimate the decay of $\ueps$ through the decay character. In Remark \ref{remark-comparison} we compare these estimates with similar ones obtained for the Navier-Stokes equations by Bjorland and M.E. Schonbek \cite{MR2493562}.  Finally, in Theorem \ref{decay-hs-rusin} we prove results concerning decay of $\dot{H}^s$ norm of solutions $\ueps$ to (\ref{eqn:rusin}). 

\section{Decay character and characterization of decay of linear systems}

\subsection{Decay Character} \label{section-decay-character} We now introduce  the main definitions used to understand the behavior of $\widehat{\Lambda^s u}_0 (\xi)$ at the origin, through comparison with $f(\xi) = |\xi|^{r}$.

\bd
Let  $u_0 \in L^2(\R^n)$, with $s \geq 0$ and let $\Lambda = (- \Delta) ^{\frac{1}{2}}$. The {\em decay indicator}  $P_r^s(u_0)$ corresponding to $\Lambda^s u_0$ is defined by

\begin{displaymath}
P_r^s(u_0) = \lim _{\rho \to 0} \rho ^{-2r-n} \int _{B(\rho)} \large| \xi|^{2s}|\widehat{u}_0 (\xi) \large|^2 \, d \xi
\end{displaymath}
for $r \in \left(- \frac{n}{2} +s, \infty \right)$, where $B(\rho)$ is the ball at the origin with radius $\rho$.
\ede

\br \, Note that $P_r^s( u_0) = 0$ for $r \leq s - \frac{n}{2}$. By setting $r = q + s$, we see that the decay indicator compares $|\widehat{\Lambda^s u}_0 (\xi)|^2$ to $f(\xi) = |\xi|^{2(q + s)}$ near $\xi = 0$. When $s = 0$, we recover the definition from Bjorland and M. E. Schonbek \cite{MR2493562}. 
\er

\bd  \label{df-decay-character} The {\em decay character of $\Lambda^s u_0$}, denoted by $r_s^{\ast} = r_s^{\ast}( u_0)$ is the unique  $r \in \left( -\frac{n}{2}+s, \infty \right)$ such that $0 < P_r^s(u_0) < \infty$, provided that this number exists. If such  $P_r^s( u_0)$ does not exist, we set $r_s^{\ast} = - \frac{n}{2}+s$, when $P_r^s(u_0)  = \infty$ for all $r \in \left( - \frac{n}{2}+s, \infty \right)$  or $r_s^{\ast} = \infty$, if $P_r^s(u_0)  = 0$ for all $r \in \left( -\frac{n}{2}+s, \infty \right)$. 
\ede

\br \, When $s=0$ we denote $P_r^0( u_0) =P_r( u_0)$, and $r_0^{\ast} = r^{\ast} $. \er

\bex \, Let $u_0 \in L^2 (\R^n)$ such that $\widehat{u}_0 (\xi) = 0$, for $|\xi| < \delta$, for some $\delta > 0$. Then, $P_r ^s (u_0) = 0$, for any $r \in \left( - \frac{n}{2} + s, \infty\right)$ and  $r^{\ast} _s (u_0) = \infty$.
\eex

\bex Let $u_0 \in L^p (\R^n) \cap L^2 (\R^n)$, with $1  \leq p \leq 2$. Then, as $\widehat{u}_0 \in L^q (\R^n)$, with $\frac{1}{p} + \frac{1}{q} = 1$, we have that 

\begin{displaymath}
\int _{B(\rho)} \large| \widehat{u}_0 (\xi) \large|^2 \, d \xi \leq \left( \int _{B(\rho)} \large| \widehat{u}_0 (\xi) \large|^{q} \, d \xi \right) ^{\frac{2}{q}} \left(  \int _{B(\rho)} \, d \xi \right) ^{2 - \frac{2}{q}}  \leq C \rho ^{2n \left( 1 - \frac{1}{p} \right)}.
\end{displaymath}
From the definition of decay character we obtain that $r^{\ast} (u_0) = - n \left( 1 - \frac{1}{p} \right)$. So, if $u_0 \in L^1 (\R^n) \cap L^2 (\R^n)$ we have that $r^*(u_0) =0$ and if $u_0 \in L^2 (\R^n)$ but $u_0 \notin L^p (\R^n)$, for any $1 \leq p < 2$, we have that $r^*(u_0) = -  \frac{n}{2}$.
\eex

\subsection{Linear part}
\label{section-linear-part}
Let $X$ be a Hilbert space on $\R^n$. Consider a pseudodifferential operator $\mathcal{L}: X^n \to \left( L^2 (\R^n) \right) ^n$, with symbol $ \Mm(\xi)$ such that 

\be
\label{eqn:symbol}
\Mm(\xi) = P^{-1} (\xi) D(\xi) P(\xi), \qquad \xi-a.e.
\ee
where $P(\xi) \in O(n)$ and $D(\xi) = - c_i |\xi|^{2\a} \delta _{ij}$, for $c_i > c>0$ and $0 < \a \leq 1$.

Given the linear equation

\be
\label{eqn:linear-part}
v_t = \mathcal{L} v.
\ee
it follows that

\begin{eqnarray*}
\frac{1}{2} \frac{d}{dt} \Vert v(t) \Vert _{L^2} ^2 & = &  \langle \widehat{v}, \Mm \widehat{ v} \rangle _{L^2} =  \langle \widehat{v}, P^{-1}  D P \widehat{v} \rangle _{L^2} \nonumber \\ & = & - \langle (-D) ^{\frac{1}{2}} P \widehat{v} ,  (-D) ^{\frac{1}{2}} P \widehat{v} \rangle _{L^2} \nonumber \\ & = & - \int _{\R^n} |(-D) ^{\frac{1}{2}} P \widehat{v}|^2 \, d \xi \nonumber \\ & \leq & - C  \int _{\R^n} |\xi|^{2 \a} |\widehat{v}|^2 \, d \xi.
\end{eqnarray*}
This is (\ref{eqn:ineq-fourier-splitting}), which we need for using the Fourier Splitting method.

\br \, The method also works, with appropiate minor changes, for $D(\xi) = - c_i |\xi|^{2\a _i} \delta _{ij}$, for $c_i \geq 0$ and $0 \leq \a_i \leq 1$, with at least one pair $(c_i, \a_i) \neq (0,0)$. 
\er

\bex \label{example-fractional-laplacian} Let $\mathcal{L}$ be the fractional Laplacian acting on vector fields on $\R^n$ by

\begin{displaymath}
\left( \mathcal{L} u \right) _i = (- \Delta) ^{\a} u_i, \qquad i = 1, \cdots n.
\end{displaymath}
Its symbol $\left( \Mm(\xi) \right) _{ij}= - C |\xi|^{2 \a} \delta _{ij}$ verifies the required condition with $P(\xi) = Id$ and $D(\xi) = - |\xi|^{2\a} \delta _{ij}$.
\eex

\bex \label{example-rusin}  Let 

\begin{displaymath}
u_t   = \mathcal{L} u = \Delta u + \frac{1}{\eps} \, \nabla div \, u , \qquad \eps > 0
\end{displaymath}
be the compressible approximation to the Stokes system in $\R^3$ introduced by Temam \cite{MR0237972}. The symbol for this operator  is $ \left( \Mm(\xi) \right) _{ij} = - |\xi|^2 \delta_{ij}  - \frac{1}{\eps} \xi_i \xi _j $, with $D(\xi) = diag (- |\xi|^2, - |\xi|^2, - \left( 1 + \frac{1}{\eps} \right) |\xi|^2 )$ and

\begin{displaymath}
P(\xi) =  \left( \begin{array} {ccc}
\frac{- \xi_2}{\sqrt{\xi_1 ^2 + \xi_2 ^2}} & \frac{- \xi_1 \xi_3}{\sqrt{1 - \xi_3 ^2}} & \xi_1 \\ \frac{\xi_1}{\sqrt{\xi_1 ^2 + \xi_2 ^2}} &  \frac{- \xi_2 \xi_3}{\sqrt{1 - \xi_3 ^2}} & \xi_2 \\ 0 &  \frac{1 - \xi_3 ^2}{\sqrt{1 - \xi_3 ^2}} & \xi_3
\end{array} \right),
\end{displaymath}
where $v = (\xi_1, \xi_2, \xi_3)$ has norm one. Then

 \be
\label{eqn:sol-fund-frec}
\left(  e^{t \Mm (\xi)} \right) _{ij} = e ^{-t |\xi|^2} \delta _{ij} - \frac{\xi _i \xi _j}{|\xi|^2} \left(e ^{-t |\xi|^2} - e ^{- \left( 1 + \frac{1}{\eps} \right)t |\xi|^2} \right),
\ee
see Rusin \cite{rusin}.
\eex

\subsection{  $L^2$ decay characterization of solutions to  linear systems (\ref{eqn:linear-part})}
\label{section-decay-l2-linear-part}

\bt
\label{characterization-decay-l2}
Let $v_0 \in L^2 (\R^n)$ have decay character $r^{\ast} (v_0) = r^{\ast}$. Let $v (t)$ be a solution to  (\ref{eqn:linear-part}) with data $v_0$. Then:
\begin{enumerate}
\item if $- \frac{n}{2 } < r^{\ast}< \infty$, there exist constants $C_1, C_2> 0$ such that
\begin{displaymath}
C_1 (1 + t)^{- \frac{1}{\a} \left( \frac{n}{2} + r^{\ast} \right)} \leq \Vert v(t) \Vert _{L^2} ^2 \leq C_2 (1 + t)^{- \frac{1}{\a} \left( \frac{n}{2} + r^{\ast} \right)};
\end{displaymath}
\item if $ r^{\ast}= - \frac{n}{2}$, there exists $C = C(\eps) > 0$ such that
\begin{displaymath}
\Vert v(t) \Vert _{L^2} ^2 \geq C (1 + t)^{-\eps}, \qquad \forall \eps > 0,
\end{displaymath}
i.e. the decay of $\Vert v(t) \Vert _{L^2} ^2$ is slower than any uniform  algebraic rate;
\item if $r^{\ast} = \infty$, there exists $C > 0$ such that
\begin{displaymath}
\Vert v(t) \Vert _{L^2} ^2 \leq C (1 + t) ^{- m}, \qquad \forall m > 0,
\end{displaymath}
i.e. the decay of $\Vert v(t) \Vert _{L^2}$ is faster than any algebraic rate.
\end{enumerate} \et
{\bf Proof:} \, (1) Suppose $P_r (u_0) > 0$. Then there exists $  \rho_0>0, C_1 > 0$ such that for $0 < \rho \leq \rho_0$ we have

\begin{displaymath}
C_1 < \rho^{-2r - n} \int _{B(\rho)} |\widehat{v_0} (\xi|^2 \, d \xi.
\end{displaymath}
Let $B( \rho(t)) = \{\xi \in \R^n:|\xi| \leq \rho (t) \}$, for some nonincresing, continuous radius $\rho = \rho (t)$ to be determined later. From (\ref{eqn:symbol}) we obtain 

\begin{displaymath}
|e^{{\mathcal M}(\xi)t} \widehat{v_0} (\xi)| = |P^{-1} (\xi) e^{D(\xi) t} P(\xi) \widehat{v_0} (\xi)| \geq C e^{-c |\xi|^{2\a}t} |\widehat{v_0} (\xi)| \geq C e^{-ct \rho^{2 \a} (t)} |\widehat{v_0} (\xi)|.
\end{displaymath}
Then 
\begin{align} \label{2.9.1}
\Vert v(t) \Vert _{L^2} ^2  \geq  \int _{B(\rho (t))}  |e^{{\mathcal M}(\xi)t} \widehat{v_0} (\xi)|^2 \, d \xi \\
\geq  C \rho ^{2r + n} e^{-ct \rho^{2 \a} (t)} \rho ^{-2r -n} \int _{B(\rho (t))} |\widehat{v_0} (\xi)|^2 \, d \xi  \geq  C \rho^{2r + n} e^{-t \rho^{2 \a} (t)}.\notag
\end{align}
Taking $\rho(t) = \rho_0 (1 + t) ^{- \frac{1}{2 \a}}$ yields $e^{-ct\rho^{2 \a} (t)} \geq C > 0$, so

\begin{displaymath}
\Vert v(t) \Vert _{L^2} ^2 \geq C (1 + t)^{- \frac{1}{\a} \left( \frac{n}{2} + r^{\ast} \right)}
\end{displaymath}
which leads to lower bound we wanted to prove.
 The upper bound follows by Fourier splitting. From (\ref{eqn:ineq-fourier-splitting}) we have

\begin{equation} \label{2.9.2}
\frac{d}{dt} \Vert v(t) \Vert _{L^2} ^2 \leq - C \int _{\R^n} |\xi|^{2 \a} |\widehat{v}(\xi)|^2  \, d \xi \leq - C \rho^{2 \a} (t) \int _{B ^c (\rho (t))} |\widehat{v} (\xi)|^2 \, d \xi,
\end{equation}
with $B(\rho(t))$ as before. Hence

\be
\label{eqn:standard-fourier-splitting}
\frac{d}{dt} \Vert v(t) \Vert _{L^2} ^2 + \rho^{2 \a} (t) \Vert v(t) \Vert _{L^2} ^2 \leq C \rho^{2 \a} (t) \int _{B(t)} |\widehat{v} (\xi)|^2 \, d \xi.
\ee
As $P_r (u_0) < \infty$, there exist $\rho_0 > 0, C > 0$ such that for $0 < \rho < \rho_0$

\be
\label{eqn:inequality-decay-indicator}
\rho^{-2r - n} \int _{B(\rho(t))} |\widehat{v_0} (\xi)|^2 \, d \xi \leq C.
\ee
Also, we have that

\be
\label{eqn:add-substract}
\int _{B(\rho (t))} |\widehat{v} (\xi)|^2 \, d \xi \leq  C \int _{B(\rho (t))} |\widehat{v_0} (\xi)|^2 \, d \xi.
\ee
Then from (\ref{eqn:standard-fourier-splitting}), (\ref{eqn:inequality-decay-indicator}) and (\ref{eqn:add-substract}) we obtain

\begin{displaymath}
\frac{d}{dt} \Vert v(t) \Vert _{L^2} ^2 \rho^{2 \a} (t) \leq C \rho^{2 \a + 2r + n} (t).
\end{displaymath}
We choose $\rho(t) = m^{\frac{1}{2 \a}} (1 + t)^{- \frac{1}{2 \a}}$, with $m > r + \frac{n}{2}$ and multiply this inequality by the integrating factor $h(t) = (1 + t)^m$ to obtain
\begin{displaymath}
\frac{d}{dt} \left( (1 + t)^m \Vert v(t) \Vert _{L^2} ^2 \right) \leq C (1 + t) ^{m - 1 - \frac{r}{\a} - \frac{n}{2 \a}}.
\end{displaymath}
Integrating between $0$ and $t$  we obtain the upper bound.

\medskip
\noindent (2) If $r^{\ast} = - \frac{n}{2}$, for any fixed $r \in (- \frac{n}{2}, \infty)$ we have $P_r (u_0) = \infty$. Then for any $\widetilde{C} = \widetilde{C} (r)$ there exists $\rho_0 > 0$ such that for $0 < \rho_0 < \rho$

\begin{displaymath}
\widetilde{C} < \rho ^{-2r - n} \int _{B(\rho)} |\widehat{u}_0 (\xi)|^2 \, d \xi.
\end{displaymath}
Proceeding as in the proof of the lower bound in (1) with an inequality similar to (\ref{2.9.1}), we obtain

\begin{displaymath}
C (1 + t)^{-r - \frac{n}{2}} \leq \Vert v(t) \Vert _{L^2} ^2.
\end{displaymath}
As this holds for any $r \in (- \frac{n}{2}, \infty)$, the estimate is established.

\noindent (3) As $r^{\ast} = \infty$, for any fixed $r \in (- \frac{n}{2}, \infty)$ we have $P_r (u_0) = 0$. Then for any $\widetilde{C} = \widetilde{C} (r)$ there exists $\rho_0 > 0$ such that for $0 < \rho_0 < \rho$

\begin{displaymath}
\rho ^{-2r - n} \int _{B(\rho)} |\widehat{u}_0 (\xi)|^2 \, d \xi < \widetilde{C}
\end{displaymath}
Proceeding as in the proof of the upper bound in (1)   with an inequality similar to (\ref{2.9.2}), we obtain

\begin{displaymath}
\Vert v(t) \Vert _{L^2} ^2 \leq C (1 + t)^{-r - \frac{n}{2}} .
\end{displaymath}
As this holds for any $r \in (- \frac{n}{2}, \infty)$, the estimate is established. $\Box$

\subsection{Characterization of decay in $\dot{H}^s$ of solutions to (\ref{eqn:linear-part})} \label{section-decay-hs-solution} In the next Theorem we establish a relation between the decay character of $\Lambda^s u_0$ and that of $u_0$.

\bt
\label{decay-character-hs}
Let $u_0 \in H^s (\R^n), s > 0$.  
\begin{enumerate}
 \item If $-\frac{n}{2} < r^{\ast} (u_0) < \infty$ then  $- \frac{n}{2} +s< r_s^{\ast}(u_0) < \infty$ and  $r_s^{\ast}(u_0) = s + r^{\ast} (u_0)$. 
\item $r_s^{\ast}  (u_0) = \infty$ if an only if $r^{\ast} (u_0) = \infty$.
\item  $r^{\ast} (u_0) =- \frac{n}{2}$ if and only if $r_s^{\ast}(u_0) = r^{\ast} (u_0) + s = - \frac{n}{2} + s$. \end{enumerate}
\et

{\bf Proof:} (1) Let $r = q +s$, for some $q>-\frac{n}{2}$.   We have that

\begin{displaymath}
\rho ^{-2r - n} \int _{B(\rho)} |\xi|^{2s}|\widehat{u}_0|^2 \, d \xi \leq \rho ^{- 2q - n} \int _{B(\rho)} |\widehat{u}_0|^2 \, d \xi,
\end{displaymath}
which,   after letting $\rho$ go to zero, leads to

\begin{displaymath}
0 \leq   P_r^s(u_0) \leq P_q (u_0) < \infty.
\end{displaymath}
If $P_r^s(u_0) > 0$, we are done, as $r_s^{\ast}  (u_0) = r = q + s = r^{\ast}  ( u_0) + s$. If $ P_r ^s (u_0)=0$, there are two possible scenarios:
\begin{enumerate}
\item[(a)] there exists $R > r $ so that $ P_R^s(u_0) >0$, or
\item[(b)]  $P_R^s(u_0) =0,\;\forall\; R\geq r.$
\end{enumerate} 
If $(a)$ is true, we are done by choosing $R = s + q$, since then
\[ 0 <  P_R^s(u_0) \leq P_q (u_0) < \infty,\]
and we can proceed as  in the last situation. 

We show now that $(b)$ leads to a contradiction.  From (b) and the decay character definition, since $0\leq P_m^ s \leq P_R ^ s,\; \forall m\leq R$,  it follows that $P_m^s ( u_0) = 0$, for all $m \in \R$. Let  $r=q+ s$. Then, $P_r^s ( u_0) = 0$ implies that given any $\e >0 $,  there is $\rho_0$ so that for all $\rho \leq \rho_0$ we have

\begin{displaymath}
\rho^{-2r-n}\int _{B(\rho)} |\xi|^{2s} |\widehat{u}_0 (\xi)|^2 \, d \xi <\e.
\end{displaymath}
Let 
\begin{displaymath}
A_p(\rho) = \{ \xi: \frac{\rho}{2^{p+1}} \leq |\xi| \leq \frac{\rho}{2^p} \}.
\end{displaymath}
Note that $A_p(\rho) \subset B(\frac{\rho}{2^p})$. Hence  by choosing  $\rho \leq \rho_0$, since $\frac{\rho}{2^p}\leq \rho\leq \rho_0$ 
\begin{eqnarray}\label{anillo1}
\left(\frac{\rho}{2^p}\right)^{-2r-n}\left(\frac{\rho}{2^{p+1}} \right)^{2s}      \int _{A_p(\rho)} |\widehat{u}_0(\xi)|^2 \, d \xi\leq \\
\left(\frac{\rho}{2^p}\right)^{-2r-n}\int _{A_p(\rho)} |\xi|^{2s} |\widehat{u}_0 (\xi)|^2 \, d \xi \leq \notag\\
\left(\frac{\rho}{2^p}\right)^{-2r-n}\int _{B(\frac{\rho}{2^p})} |\xi|^{2s} |\widehat{u}_0 (\xi)|^2 \, d \xi <\e.\notag
\end{eqnarray}
Since $r= q+s$ the first term from the  above  inequalities can be expressed as
\[\left(\frac{\rho}{2^p}\right)^{-2r-n}\left(\frac{\rho}{2^{p+1}} \right)^{2s}      \int _{A_p(\rho)} |\widehat{u}_0(\xi)|^2 \, d \xi=
 2^{p(2q+n)} 2^{-2s} \rho^{-2q-n}   \int _{A_p(\rho)} |\widehat{u}_0(\xi)|^2 \, d \xi\]
Combining the last equality with (\ref{anillo1}) yields

\[ \rho^{-2q-n}  \int _{A_p(\rho)} |\widehat{u}_0(\xi)|^2 \, d \xi \leq \eps \, 2^{2s} 2^{-p(2q+n)}. \]
Note that $\cup_{p=0}^{\infty} A_p(\rho)=B(\rho)$. Summing over $p$ yields

\[ \rho^{-2q-n}  \int _{B(\rho)} |\widehat{u}_0(\xi)|^2 \, d \xi \leq \eps \, 2^{2s} \sum_ {p=0}^{\infty}2^{-p(2q+n)}. \]
Hence 

\[ \rho^{-2q-n}  \int _{B(\rho)} |\widehat{u}_0(\xi)|^2 \, d \xi \leq \eps \, 2^{2s}\frac{1}{1 - 2^{2q+n}}. \]
Since $\e>0$ was arbitrary it follows that
\[ P_r(u_0) = \lim_{\rho \to 0} \rho^{-2q-n}  \int _{B(\rho)} |\widehat{u}_0(\xi)|^2 =0\]
Our hypothesis was that $r^*> -\frac{n}{2}$, hence by definition  $P_r(u_0)>0$ and we reached a contradiction.  This completes the proof of part (1).

\smallskip

\noindent (2)  Let $r \in \left( - \frac{n}{2} + s, \infty \right)$. Then

\begin{displaymath}
\rho^{-2r - n} \int _{B(\rho)} |\xi|^{2s} |\widehat{u}_0 (\xi)|^2 \, d \xi \leq \rho^{-2(r -s) - n} \int _{B(\rho)} |\widehat{u}_0 (\xi)|^2 \, d \xi.
\end{displaymath}
Taking limits when $\rho$ goes to zero we obtain $0 \leq P_r ^s (u_0) \leq P_{r - s} (u_0)$. If $r^{\ast} (u_0) = \infty$, then $P_{r - s} (u_0) = 0$, for any $r - s \in \left( - \frac{n}{2}, \infty \right)$. So $P_r ^s (u_0) = 0$, hence $r^{\ast} _s (u_0) = \infty$.

Now suppose $r^{\ast} _s (u_0) = \infty$, but $r^{\ast} (u_0) = m < \infty$. Then, by part (1) we have that $r^{\ast} _s (u_0) = m + s < \infty$, which is a contradiction.

\smallskip

\noindent (3)  Let $r \in \left( - \frac{n}{2} + s, \infty \right)$. Then

\begin{displaymath}
\rho^{-2r - n} \int _{B(\rho)} |\xi|^{2s} |\widehat{u}_0 (\xi)|^2 \, d \xi \leq \rho^{-2(r -s) - n} \int _{B(\rho)} |\widehat{u}_0 (\xi)|^2 \, d \xi.
\end{displaymath}
Suppose $r^{\ast} _s (u_0) = - \frac{n}{2} + s$. Then, taking limits when $\rho$ goes to zero we obtain

\begin{displaymath}
\infty = P_r ^s (u_0) \leq P_{r - s} (u_0).
\end{displaymath}
As this holds for any $r - s \in \left( - \frac{n}{2}, \infty \right)$, we conclude that $r^{\ast} (u_0) = - \frac{n}{2}$.

Now suppose $r^{\ast} (u_0) = - \frac{n}{2}$ and suppose that $r^{\ast} _s (u_0) = q >  - \frac{n}{2} + s$. If $q = \infty$, then by part (2) we would have $r^{\ast} (u_0) = \infty$, which is a contradiction. If $q < \infty$, then from

\begin{displaymath}
\rho^{-2r - n} \int _{B(\rho)} |\xi|^{2s} |\widehat{u}_0 (\xi)|^2 \, d \xi \leq \rho^{-2(r -s) - n} \int _{B(\rho)} |\widehat{u}_0 (\xi)|^2 \, d \xi.
\end{displaymath}
we obtain $0 < P_q ^s (u_0) \leq P_{q - s} (u_0)$. Then $r^{\ast} (u_0) = q - s > - \frac{n}{2}$, which is also a contradiction. Then we must have $r^{\ast} _s (u_0) = - \frac{n}{2} + s$. $\Box$

This Theorem  leads to the following result.

\bt
\label{characterization-decay-hs}
Let $v_0 \in H^s(\R^n), s > 0$ have decay character $r^{\ast} _s = r^{\ast} _s (v_0)$. Then

\begin{enumerate}
\item if $- \frac{n}{2 } \leq r^{\ast} < \infty$, there exist constants $C_1, C_2 > 0$ such that

\begin{displaymath}
C_1 (1 + t)^{- \frac{1}{\a} \left( \frac{n}{2} + r^{\ast}  + s \right)} \leq \Vert v(t) \Vert _{\dot{H}^s} ^2 \leq C_2 (1 + t)^{- \frac{1}{\a} \left( \frac{n}{2} + r^{\ast}  + s \right)};
\end{displaymath}
\item if $r^{\ast} = \infty$, then

\begin{displaymath}
\Vert  v(t) \Vert _{\dot{H}^s} ^2 \leq C (1 + t) ^{- r}, \qquad \forall r > 0,
\end{displaymath}
i.e. the decay of $\Vert  v(t) \Vert _{\dot{H}^s}$ is faster than any algebraic rate.
\end{enumerate}
\et

{\bf Proof:} \, From $(2)$ and $(3)$ in Theorem \ref{decay-character-hs}, we have that if $- \frac{n}{2 } \leq r^{\ast} < \infty$ then $r^{\ast} (v_0) = s + r^{\ast}$. We then can apply $(1)$ in  Theorem \ref{characterization-decay-l2} with initial data $\Lambda ^s v_0 \in L^2 (\R^n)$. If $r^{\ast} = \infty$, we use $(1)$ in Theorem \ref{decay-character-hs} and $(3)$ in  Theorem \ref{characterization-decay-l2}. $\Box$

\section{Applications: Dissipative Quasi-Geostrophc Equation and Compressible Approximation to Navier-Stokes Equations}
\label{applications}

\subsection{Dissipative Quasi-Geostrophic Equation}
\label{qg}

Consider the dissipative quasi-geostrophic equation  in $\R^2$

\be
\label{eqn:qg}
\theta _t + u \cdot \nabla \theta + \kappa (- \Delta) ^{\a}u = 0, \qquad 0 < \a \leq 1, \kappa > 0
\ee
where $u = R^{\bot} \theta = \left( -R_2 \theta, R_1 \theta \right)$, $R_i$ is the Riesz transform in the $x_i$ variable and $\theta = \theta (x,t)$ is a scalar, the potential temperature of the fluid with velocity $u$. This active scalar equation has been extensively studied in the last 20 years, since for $\a = \frac{1}{2}$, it provides a dimensionally correct $2D$ model to the $3D$ Navier-Stokes equations. Moreover, this equation and its inviscid counterpart (i.e., $\kappa = 0$ in (\ref{eqn:qg})), describe important models of currents and atmosphere circulation, see Majda and Tabak \cite{MR1422288}, Pedlosky \cite{Pedlosky}. An extensive bibliography concerning  both the inviscid and dissipative equations arising from (\ref{eqn:qg}), can be found in Chae, Constantin, D. C\'ordoba, Gancedo and Wu \cite{CPA:CPA21390} and Chae, Constantin and Wu \cite{chae-constantin-wu-2}, \cite{MR3119608}.

The decay of the $L^2$ and Sobolev norms and asymptotic behaviour of solutions to this equation has been addressed in many articles, see for example Benameur and Blel \cite{MR2922933}, Carrillo and Ferreira \cite{MR2295181}, \cite{MR2322939}, \cite{MR2412324}; Constantin and Wu \cite{MR1709781}; A. C\'ordoba and D. C\'ordoba \cite{MR2084005}; Dong and Du \cite{MR2399451}; Niche and Planas \cite{MR2914582}; Niche and M. E. Schonbek \cite{MR2342289}; M. E. Schonbek and T. Schonbek \cite{MR2001105}, \cite{MR2166670}; Tun and Miyakawa \cite{MR2730620} and Zhou \cite{MR2171725}, \cite{MR2430660}. As its linear operator has the properties needed to use the results obtained in Section \ref{section-decay-character} (see Example \ref{example-fractional-laplacian}), we establish decay estimates for  (\ref{eqn:qg}) in terms of the decay character of the initial data $\theta_0 \in L^2$. 

\subsubsection{Decay of the $L^2$ norm of solutions to the Quasi-Geostrophic Equation}

\bt
\label{decay-l2-qg}
Let $\theta_0 \in L^2 (\R^2)$, let $r^{\ast} = r^{\ast} (\theta_0), -1 < {r^{\ast}} < \infty$, and $0< \a\leq 1$.  Let $\theta$ be a weak solution to (\ref{eqn:qg}) with data $\theta_0$. Then:
\begin{enumerate}
\item If $r^* \leq 1-\a$, then 
\begin{displaymath}
\|\theta(t)\|_{L^2} ^2 \leq C(t+1)^{-\frac{1}{\a} (1 +r^{\ast})};
\end{displaymath}
\item if $r^* \geq 1-\a$, then 
\begin{displaymath}
\|\theta(t)\|_{L^2} ^2 \leq C(t+1)^{- \frac{1}{\a}\left(2 - \a \right)}.
\end{displaymath}
\end{enumerate}
\et

{\bf Proof:} Existence of weak solutions to (\ref{eqn:qg}) was established by Resnick \cite{Resnick}. When using the Fourier Splitting method we  prove the estimates for the full nonlinear equations, assuming the solutions are regular enough. For full details of the limiting argument used to obtain the estimate for to weak solutions, see pages 267--269 in Lemari\'e-Rieusset \cite{MR1938147} and the Appendix in Wiegner \cite{MR881519}.  

Let

\begin{displaymath}
B(t) = \{\xi \in \R^2: |\xi|^{2\a} \leq \frac{f'(t)}{2 f(t)} \}.
\end{displaymath}
The Fourier Splitting method, yields 
\begin{align*}
\frac{d}{dt} \left( f(t) \Vert \theta (t) \Vert _{L^2} ^2 \right) \leq f'(t) \int _{B(t)} |\widehat{u} (\xi, t)|^2 \, d \xi \\\leq C f'(t) \left( \Vert \Theta (t) \Vert _{L^2} ^ 2 + \int _{B(t)}  \left( \int _0 ^t e^{- (t - s) |\xi|^{2 \a}} |\xi| |\widehat{u \theta} (\xi, s)|^2 \, ds \right) ^2 \, d \xi \right)
\end{align*}
where $\Theta$ is the solution to the linear part. We have the straightforward inequality

\begin{displaymath}
\int _0 ^t e^{- (t - s) |\xi|^{2 \a}} |\xi||\widehat{u \theta} (\xi, s)|^2 \, ds \leq \int _0 ^t |\xi| \Vert \theta (s) \Vert _{L^2} ^2 \, ds
\end{displaymath}
which leads to

\begin{displaymath}
\left( \int _0 ^t e^{- (t - s) |\xi|^{2 \a}} |\xi||\widehat{u \theta} (\xi, s)|^2 \, ds \right) ^2 \leq t \int _0 ^t |\xi| ^2 \Vert \theta (s) \Vert _{L^2} ^4 \, ds.
\end{displaymath}
Thus
\begin{displaymath}
\int _{B(t)} \left( \int _0 ^t e^{- (t - s) |\xi|^{2 \a}} |\xi||\widehat{u \theta} (\xi, s)|^2 \, ds \right) ^2 \, d \xi \leq C t \left( \frac{f'(t)}{f(t)} \right) ^{\frac{2}{\a}} \int _0 ^t\Vert \theta (s) \Vert _{L^2} ^4 \, ds, 
\end{displaymath}
hence
\be
\label{eqn:main-inequality-geostrophic}
\frac{d}{dt} \left( f(t) \Vert \theta (t) \Vert _{L^2} ^2 \right) \leq f'(t) \left( \Vert \Theta (t) \Vert _{L^2} ^ 2 + t \left( \frac{f'(t)}{f(t)} \right) ^{\frac{2}{\a}} \int _0 ^t\Vert \theta (s) \Vert _{L^2} ^4 \, ds \right).
\ee
We first obtain a preliminary decay, which will be later used to obtain the optimal decay rate. Consider $0 < \a < 1$ and  let $f(t) = [\ln (e + t)] ^{1 + \frac{1}{\a}}$. Then

\begin{displaymath}
t f'(t)\left( \frac{f'(t)}{f(t)} \right) ^{\frac{2}{\a}} \int _0 ^t\Vert \theta (s) \Vert _{L^2} ^4 \, ds \leq C \frac{(e + t)}{[\ln (e + t)]^{\frac{1}{\a}}(e + t)^{\frac{2}{\a}}}\leq C \frac{1}{(e+t)^{\frac{2}{\a}-1}}.
\end{displaymath}
which is integrable since $0 < \a < 1$. Integration of (\ref{eqn:main-inequality-geostrophic}) combined with Theorem \ref{characterization-decay-l2} yields

\begin{eqnarray*}
\Vert \theta (t) \Vert _{L^2} ^2 & \leq &  \Vert \Theta (t) \Vert _{L^2} ^2 + C [\ln (e + t)] ^{- \left( 1 + \frac{1}{\a} \right)} \nonumber \\ & \leq & C (t + 1) ^{- \left( \frac{1 + r^{\ast}}{\a} \right)} +  C [\ln (e + t)] ^{- \left( 1 + \frac{1}{\a} \right)} \leq C [\ln (e + t)] ^{- \left( 1 + \frac{1}{\a} \right)}.
\end{eqnarray*}
Now for $\a = 1$, take $f(t) = [\ln (e + t)]^3$. Then

\begin{displaymath}
t \left( \frac{f'(t)}{f(t)} \right) ^{\frac{2}{\a}} \int _0 ^t\Vert \theta (s) \Vert _{L^2} ^4 \, ds \leq C [\ln (e + t)] ^{-2}
\end{displaymath}
which, by the same argument as before, leads to the same decay as before, in this case with $\a = 1$, 

\be\label{ldecay}
\Vert \theta (t) \Vert _{L^2} ^2 \leq C [\ln (e + t)] ^{-2}.
\ee
We  now proceed by bootstrapping. Let $0 < \a \leq 1$. From the preliminary  decay we have

\begin{displaymath}
\Vert \theta (t) \Vert _{L^2} ^4 \leq \frac{ C\Vert \theta(t) \Vert _{L^2} ^2}{[\ln (e + t)] ^{1 + \frac{1}{\a}}}.
\end{displaymath}
Let $\b > 1+  \frac{1}{\a} + \frac{r^{\ast}}{\a} $. Choose $f(t) = (t + 1) ^{\b}$. Then, plugging in (\ref{eqn:main-inequality-geostrophic}) and integrating in time, we obtain

\begin{eqnarray*}
(t + 1)^{\b} \Vert \theta (t) \Vert _{L^2} ^2 & \leq & \Vert \theta_0 \Vert _{L^2} ^2 + \int _0 ^t (s + 1)^{\b- 1} \Vert \Theta (s) \Vert _{L^2} ^2 \, ds \\ & + & \int _0 ^t (s + 1)^{\b- 1} s(s+ 1) ^{- \frac{2}{\a}} \int _0 ^t \Vert \theta (\tau) \Vert _{L^2} ^4 \, d \tau ds\\  & \leq & \Vert \theta_0 \Vert _{L^2} ^2 + \int _0 ^t (s + 1)^{\b - 1 - \frac{1}{\a} - \frac{r^{\ast}}{\a}} \, ds + (t+1)^{\b -\frac{2}{\alpha} +1}\int _0 ^t \Vert \theta (\tau) \Vert _{L^2} ^4 \, d \tau  \end{eqnarray*}
where we used Theorem \ref{characterization-decay-l2} to estimate the decay of the linear part $\Theta (t).$ Dividing by $(t+1)^{\b -\frac{2}{\alpha} +1}$  and using estimate (\ref{ldecay}) yields 
\begin{eqnarray*}
(t + 1)^{\frac{2}{\a}-1} \Vert \theta (t) \Vert _{L^2} ^2  & \leq & \Vert \theta_0 \Vert _{L^2} ^2 ( t+1)^{-(\b -\frac{2}{\alpha} +1)}   + C (t + 1) ^{\frac{1}{\a}- \frac{r^{\ast}}{\a}-1} \nonumber \\ & + & C \int _0 ^t \frac{(s + 1) ^{1- \frac{2}{\a}} }{\ln (e + s)]^{1 + \frac{1}{\a}}}(s + 1)^{\frac{2}{\a}-1}  \Vert \theta (s) \Vert _{L^2} ^2 \, ds.
\end{eqnarray*}
Since $\b$ can be taken as large as needed we also suppose that $(\b -\frac{2}{\alpha} +1)>0$. Let
\begin{align*}
\psi(t)  =  (1+t)^{\frac{2}{\a}-1} \| \theta(t) \| _{\eledos} ^2, \ a(t) = C (t+ 1) ^{- \beta} + (1 + t) ^{\frac{1}{\a}- \frac{r^{\ast}}{\a}-1}, \\ b(t)   =   C [\ln (e + t) ]^ {- \left( 1 + \frac{1}{\a} \right)}(s+1)^{1-\frac{2}{\a}}.
\end{align*}
Then the previous inequality becomes

\begin{displaymath}
\psi (t) \leq a(t) + C \int _0 ^t \psi (s) b(s) \, ds. 
\end{displaymath}
Notice that  since $0<\a\leq 1$

\begin{displaymath}
\int _0 ^t b(s) \, ds \leq C.
\end{displaymath}
We need to consider two cases:\\
{\it Case 1.} $r^* \geq 1-\a$. In this case ${\frac{1}{\a}- \frac{r^{\ast}}{\a}-1} \leq 0$. Hence $a(t) \leq C.$
Hence by a standard Gronwall inequality we have 

\begin{displaymath}
\psi (t) \leq C \, exp \left( \int _0 ^t b(s) \, ds \right)
\end{displaymath}
which yields 
\be \label{first}  \|\theta(t)\|_{L^2} ^2 \leq C(t+1)^{- \frac{1}{\a}\left(2 - \a \right)} \ed
{\it Case 2.} $r^* \leq 1-\a$. Hence  ${\frac{1}{\a}- \frac{r^{\ast}}{\a}-1} \geq 0$. In this case $a(t)$ is increasing. Thus Corollary 1.2, page 4 from Ba{\u\i}nov and Simeonov \cite{MR1171448} yields

\begin{displaymath}
\psi (t) \leq a(t) \, exp \left( \int _0 ^t b(s) \, ds \right)
\end{displaymath}
Since  
\begin{displaymath}
a(t) (t+1) ^{1-\frac{2}{\a}} \leq C (t + 1) ^{- \beta} + C (t+1)^{-(\frac{1}{\a}+\frac{r^{\ast}}{\a})},
\end{displaymath}
it follows that
\be \label{second}  \|\theta(t)\|_{L^2} ^2 \leq C(t+1)^{-\frac{1}{\a} (1 +r^{\ast})}. \ed
The conclusion of the theorem follows by (\ref{first}) and (\ref{second})  $\Box$

We  now address  the decay of the nonlinear part $w(t) = \theta (t) - \Theta (t)$. 

\bt 
\label{decay-linear-part-qg}
Let $0 < \a \leq 1$, $\theta_0 \in  L^2(\R^2)$, $r^{\ast}=r^{\ast} (\theta_0)$. Then

\begin{enumerate}
\item if $r^{\ast} \geq 1 - \a$, then
\begin{displaymath}
\Vert \theta (t) - \Theta (t) \Vert _{L^2} ^2 \leq C (1 + t) ^{- \frac{1}{\a} \min \{3 - 2 \a, 2 \}};
\end{displaymath}

\item if $\a-1 \leq r^{\ast} \leq 1-\a$ then
\begin{displaymath}
\Vert \theta (t) - \Theta (t) \Vert _{L^2} ^2 \leq C (1 + t) ^{- \frac{1}{\a} \min \{ 2, 2 - \a + r^{\ast} \}};
\end{displaymath}
 \item  Let  $-1  \leq r^{\ast} < \a-1$ then  
\begin{displaymath}
\Vert \theta (t) - \Theta (t) \Vert _{L^2} ^2 \leq C (1 + t) ^{- \frac{1}{\a} \left( 2 - \a + r^{\ast} \right)}.
\end{displaymath}
\end{enumerate}
\et

{\bf Proof:}  The proof is similar to that of Theorem 4.3 in Constantin and Wu \cite{MR1709781}. We recall that  by Proposition 4.2 in \cite{MR1709781} we have that 

\begin{displaymath}
\|\nabla \Theta (t)\|_{\infty}  \leq C(t +1)^{-\frac{1}{\a}}.
\end{displaymath}
This inequality combined  with the results of Theorem \ref{decay-l2-qg}  yields the following  estimate of  the nonlinear term  

\[ \left| \int_{\R^2} \Theta(u\cdot \nabla \theta) dx \right|  \leq \|\nabla \Theta (t)\|_{\infty} \|\theta (t) \|_2^2 \leq C(1+t)^{-\gamma}  = h(t), \]
where 

\be
\label{eqn:gamma}
\gamma = \frac{1}{\a} (2 + r^{\ast}) \,  \text{if} \, r^{\ast} \leq 1 - \a, \qquad \gamma = \frac{1}{\a}(3 - \a), \, \text{if} \, r^{\ast} \geq 1 - \a.
\ee
To use Fourier Splitting, let $B(t) = \{\xi  \in \R^2: |\xi| \leq g(t) \}$, with $g(t)$ an  increasing and  continuous function,  to be determined below, then

\be
\label{eqn:f-s-nl}
\frac{d}{dt} \Vert w(t) \Vert _{L^2} ^2 + 2 g^{2 \a} (t) \Vert w(t) \Vert _{L^2} ^2 \leq 2 g^{2 \a} (t) \int _{B(t)} |\widehat{w} (\xi, t)|^2 \, d \xi + h(t).
\ee   
Let $r^* \geq 1-\a$. Consider first  $0 < \a < 1$. Then 

\begin{displaymath}
|\widehat{w} (\xi, t)| \leq C \int _0 ^t e^{- |\xi|^{2 \a}(t - s)} |\xi| |\widehat{\theta} \ast \widehat{\theta} (\xi, s)| \, ds \leq C |\xi| \int _0 ^t \Vert \theta (s) \Vert _{L^2} ^2 \, ds,
\end{displaymath}
which leads to

\begin{eqnarray*}
\int _{B(t)} |\widehat{w} (\xi, t)|^2 \, d \xi & \leq & C \int  _{B(t)} |\xi|^2 \left( \int _0 ^t \Vert \theta (s) \Vert _{L^2} ^2 \, ds \right)^2 \, d \xi \\ & \leq &  C \int  _{B(t)} |\xi|^2 t \left( \int _0 ^t  \Vert \theta (s) \Vert _{L^2} ^4 \, ds \right) \, d \xi \\ & \leq & C  \int  _{B(t)} |\xi|^2 \, d \xi \leq C g^4 (t),
\end{eqnarray*}
where we used the decay from Theorem \ref{decay-l2-qg} in the second inequality. After using this and the integrating factor $k(t) = \exp \left( 2 \int _0 ^t g^{2 \a} (s) \, ds \right)$ in (\ref{eqn:f-s-nl}) we obtain

\begin{displaymath}
\frac{d}{dt} \left( \exp \left( 2 \int _0 ^t g^{2 \a} (s) \, ds \right) \Vert w(t) \Vert _{L^2} ^2 \right) \leq \exp \left( 2 \int _0 ^t g^{2 \a} (s) \, ds \right) \left( C g^{2 \a + 4} (t) + h(t) \right) .
\end{displaymath}
Taking $g^{2 \a} (t) = \frac{\beta}{2 (1 + t)}$ and integrating we obtain, for large enough $\beta > 0$, 

\begin{displaymath}
\Vert w(t) \Vert _{L^2} ^2 \leq C (1 + t) ^{- \beta} + C (1 + t) ^{- \frac{2}{\a}} + C h(t) (1 + t).
\end{displaymath}
The last inequality combined with the definition of $h$ and (\ref{eqn:gamma}) yields  the result.
Now let $\a = 1$ and $g^{2} (t) = \frac{3}{2(1 + t)}$. Then for $\xi \in B(t)$, the decay obtained in Theorem \ref{decay-l2-qg} yields
\begin{displaymath}
|\widehat{w} (\xi, t)| \leq C \int _0 ^t  |\xi| \Vert \theta (s) \Vert _{L^2} ^2 \, ds \leq C \int _0 ^t (s + 1) ^{- \frac{3}{2}} \, ds \leq C,
\end{displaymath}
so (\ref{eqn:f-s-nl}) becomes

\begin{eqnarray*}
\frac{d}{dt} \Vert w(t) \Vert _{L^2} ^2  +  \frac{3}{1 + t} \Vert w(t) \Vert _{L^2} & \leq & \frac{C}{1 + t}  \int _{B(t)} |\widehat{w} (\xi, t)|^2 \, d \xi + h(t) \\ & \leq & C (1 + t) ^{-2} + C (1 + t) ^{-2}.
\end{eqnarray*}
Multiplying by $ (t + 1) ^3$ and integrating yields  the conclusion. 

\medskip

\noindent {\bf Case 2:}  Consider first $r^{\ast} \leq 1 - \a$ and $\frac{1}{\a}(1+r^*)>1$. The conclusion follows proceeding  exactly as in the case $r^{\ast} \geq 1 - \a$ and $0 < \a < 1$, with $h(t)$ according to (\ref{eqn:gamma}).  For $\frac{1}{\a}(1+r^*)= 1$,   we  proceed as  in the case $\a = 1$. The details are as follows:  let $g^{2 \a} (t) = \frac{\b}{2(1 + t)}$,  for $\xi \in B(t)$ we obtain
 
\begin{displaymath}
|\widehat{w} (\xi, t)| \leq C \int _0 ^t  |\xi| \Vert \theta (s) \Vert _{L^2} ^2 \, ds \leq C \int _0 ^t (s + 1) ^{- (\frac{1}{2 \a} + 1)} \, ds \leq C,
\end{displaymath}
which, combined with (\ref{eqn:f-s-nl}), leads to

\begin{eqnarray*}
\frac{d}{dt} \Vert w(t) \Vert _{L^2} ^2  & + &  \frac{\b}{1 + t} \Vert w(t) \Vert _{L^2}  \leq  \frac{C}{(1 + t)}  \int _{B(t)} |\widehat{w} (\xi, t)|^2 \, d \xi + h(t) \\  & \leq &   C  (1 + t) ^{-(1+ \frac{3}{2\a})} + C (1 + t) ^{- \frac{1}{\a} \left( 2 + r^{\ast} \right)} \\ & \leq & C( t+1)^{-1}[ (t+1)^{-\frac{3}{2\a}}+(t+1)^{-\frac{1}{\a}}].
\end{eqnarray*}
Multiply by $i(t) = (t + 1)^{\b}$, with $\b > 0$ large enough,  and integrate to obtain the conclusion. 

\medskip 

\noindent {\bf Case 3:} Let  $\frac{1 + r^{\ast}}{\a} < 1$. As before, with $g^{2 \a} (t) = \frac{\b}{2(1 + t)}$ and $\xi \in B(t)$, we have

\begin{displaymath}
|\widehat{w} (\xi, t)| \leq C \int _0 ^t  |\xi| \Vert \theta (s) \Vert _{L^2} ^2 \, ds \leq C|\xi| \int _0 ^t \Vert \theta (s) \Vert _{L^2} ^2 \, ds \leq C (t + 1) ^{- \frac{1}{\a} (\frac{3}{2} + r^{\ast} - \a)}.
\end{displaymath}
From (\ref{eqn:f-s-nl}) we obtain

\begin{eqnarray*}
\frac{d}{dt} \Vert w(t) \Vert _{L^2} ^2  +  \frac{\b}{1 + t} \Vert w(t) \Vert _{L^2} & \leq & \frac{2}{1 + t}  \int _{B(t)} |\widehat{w} (\xi, t)|^2 \, d \xi + h(t) \\ & \leq & C (1 + t) ^{- \frac{1}{\a} \left(\frac{9}{2} + 2r^{\ast} - \a \right)} + C (1 + t) ^{- \frac{1}{\a} \left( 2 + r^{\ast} \right)}.
\end{eqnarray*}
Multiplying by $i(t) = (t + 1)^{\b}$, for large enough $\beta > 0$ and then integrating, yields the result. $\Box$

In the next Theorem we analyze the lower bounds for rates of decay. These bounds are obtained using the reverse triangle inequality and the decays proved in Theorems \ref{decay-l2-qg} and \ref{decay-linear-part-qg}. As is usual when the linear part has fast decay, we are not able to obtain lower bounds for the solution's decay, see for example Theorem A in Miyakawa and M.E. Schonbek \cite{MR1844282} and Theorem 6.5 in Bjorland and M.E. Schonbek \cite{MR2493562}.

\bt
\label{lower-bounds-qg}
Let $0 < \a \leq 1$, $\theta_0 \in  L^2(\R^2), r^{\ast}=r^{\ast} (\theta_0)$. Then, for $0 < \a \leq \frac{1}{2}$ and $-1 < r^{\ast} \leq 1$ or $\frac{1}{2} < \a \leq 1$ and $-1 < r^{\ast} \leq 2 (1 - \a)$ we have that 

\begin{displaymath}
\Vert \theta (t) \Vert _{L^2} ^2 \geq C (1 + t) ^{-\frac{1}{\a} (1 + r^{\ast})}.
\end{displaymath}
\et

{\bf Proof:} As we only have upper bounds for the decay of the difference $\theta - \Theta$, the only instance where the reverse triangle inequality would lead to a lower bound for the decay of $\theta$ is when

\begin{displaymath}
\|\theta(t)\|^2_{L^2} \geq  \|\Theta(t) \| _{L^2} ^2- \|\theta(t)-\Theta(t)||_{L^2} ^2,
\end{displaymath}
i.e. when the decay of the linear part is slower than that of the difference. The result then follows from the estimates in Theorem \ref{characterization-decay-l2}  and Theorem  \ref{decay-linear-part-qg}.  $\Box$

Combining  the estimates from Theorems \ref{decay-l2-qg} and \ref{lower-bounds-qg}, we obtain the following.
 
\bt
Let  $\theta_0 \in L^2 (\R^2)$, with decay character $r^{\ast} = r^{\ast} (u_0)$.
\begin{enumerate}\item If $r^{\ast} \leq 1 - \a$, then there exists constants $C_1, C_2 > 0$ so that
\begin{displaymath}
C_1 (1 + t) ^{-\frac{1}{\a} (1 + r^{\ast})} \leq \Vert \theta (t) \Vert _{L^2} ^2 \leq C_2 (1 + t) ^{-\frac{1}{\a} (1 + r^{\ast})};
\end{displaymath}
\item In the region determined by  $r^{\ast} \geq 1 - \a$, $r^{\ast} \leq 1$, $r^{\ast} \leq 2(1 - \a)$ and for some $C_1, C_2 > 0$ we have

\begin{displaymath}
C_1 (1 + t) ^{-\frac{1}{\a} (1 + r^{\ast})} \leq \Vert \theta (t) \Vert _{L^2} ^2 \leq C_2 (1 + t) ^{-\frac{1}{\a} (2 - \a)};
\end{displaymath}
\item In the region determined by $r^{\ast} > 1$ and $r^{\ast} \geq 2(1 - \a)$ we have that

\begin{displaymath}
\Vert \theta (t) \Vert _{L^2} ^2 \leq C_2 (1 + t) ^{-\frac{1}{\a} (2 - \a)}.
\end{displaymath}
\end{enumerate}
\et

\subsubsection{Decay of the $\dot{H}^s$ norm} 
We recall that we have existence and regularity of solutions $\theta \in  H^s(\R^2)$  with initial data in $\theta_0 \in H^s(\R^2)$,  provided $\frac{1}{2} < \a \leq 1$, see Constantin and Wu  \cite{MR1709781}. Our goal is to prove the following
\bt 
\label{decay-derivative-qg}
Let $\frac{1}{2} < \a \leq 1$, $\a \leq s$ and $\theta_0 \in H^s (\R^2)$. For $r^{\ast} = r^{\ast} (\theta_0)$ we have that:

\begin{enumerate}
\item if $r^* \leq 1-\a$, then 
\begin{displaymath}
\|\theta(t)\|_{\dot{H}^s} ^2 \leq C(t+1)^{- \frac{1}{\a} \left( s + 1 +r^{\ast} \right)};
\end{displaymath}
\item if $r^* \geq 1-\a$, then 
\begin{displaymath}
\|\theta(t)\|_{\dot{H}^s} ^2 \leq C(t+1)^{- \frac{1}{\a} \left( s + 2 - \a \right)}.
\end{displaymath}
\end{enumerate}
\et

To do that we first need the following preliminary decay

\bt 
\label{decay-preliminary-derivative-qg}
Let $\frac{1}{2} < \a \leq 1$, $\a \leq r \leq s$ and $\theta_0 \in H^s (\R^2)$. For $r^{\ast} = r^{\ast} (\theta_0)$ we have

\[\Vert \theta(t) \Vert _{\dot{H}^r} \leq \begin{cases}
   C(t+1)^{-\frac{1}{\a}(r^{\ast} +1)},& \text{if } r^{\ast} \leq 1-\a \\
     C(t+1)^{-\frac{1}{\a}(2-\a)},& \text{if } r^{\ast} \geq 1-\a 
\end{cases}.\]
\et
\pf
The proof follows closely that of Theorem 2.4 in M.E. Schonbek and T. Schonbek \cite{MR2001105}.  Everything goes  through until (2.17) in page 362, where $\Vert \theta (t) \Vert _{L^2} ^2$ is used for the first time. Estimate (2.17) is replaced  using  the decay estimates obtained for $\|\theta(t)\|_{L^2}$  in  Theorem \ref{decay-l2-qg}

\[\frac{d}{dt} \Vert \Lambda^r \theta(t) \Vert _{L^2} ^2 + C \Vert \Lambda^r \theta (t) \Vert _{L^2} ^2 \leq C \Vert \theta (t) \Vert _{L^2} ^2 \leq \begin{cases}
   C(t+1)^{-\frac{1}{\a}(r^{\ast} +1)},& \text{if } r^{\ast} \leq 1-\a \\
     C(t+1)^{-\frac{1}{\a}(2-\a)},& \text{if } r^{\ast} \geq 1-\a 
\end{cases}.\]
To finish the proof of Theorem \ref{decay-preliminary-derivative-qg} we follow the same steps as in \cite{MR2001105}. $\Box$

\medskip

\pf ({\it  Theorem \ref{decay-derivative-qg} })
 The proof is  essentially the same as the proof of Theorem 3.2 in \cite{MR2001105}. By interpolation results that follow from Theorem \ref{decay-l2-qg} and \ref{decay-preliminary-derivative-qg} we obtain that for any $2 \leq q < \infty$

\[\Vert \theta (t) \Vert _{L^q} ^2 \leq  \begin{cases}
   C(t+1)^{-\frac{1}{\a}(r^{\ast} +1)},& \text{if } r^{\ast} \leq 1-\a \\
     C(t+1)^{-\frac{1}{\a}(2-\a)},& \text{if } r^{\ast} \geq 1-\a 
\end{cases}.\]
This  is the estimate  needed in   Theorem 3.2 for inequality (3.7), page 367. From then onwards, the proof is identical. $\Box$

\subsection{Compressible Approximation to Navier-Stokes Equations}

\label{rusin-section}

Incompressibility in the Navier-Stokes equations (\ref{eqn:navier-stokes}) leads to the nonlocal relation $p = \Delta ^{-1} \left( \partial_j v^i \partial_i v^j\right)$ between the pressure and the velocity. This nonlocality poses an important problem when trying to use numerical methods to study the behaviour of the solutions to (\ref{eqn:navier-stokes}). In order to avoid these problems, Temam \cite{MR0237972} introduced an approximation to the Navier-Stokes equations by directly relating the pressure and the velocity through  $\eps \, p = - div \, u$. Moreover, to ``stabilize'' the system, i.e. to have an energy inequality, he added the nonlinear term $\frac{1}{2} (div \, \ueps) \ueps$. This lead to the compressible system

\begin{eqnarray}  \label{eqn:rusin}
\partial _t \ueps  + (\ueps \cdot \nabla) \ueps  + \frac{1}{2} (div \, \ueps) \ueps & = & \Delta \ueps + \frac{1}{\eps} \, \nabla \cdot div \, \ueps  \nonumber \\ \ueps (x,0)  & = & \ueps _0 (x).
\end{eqnarray}
This system has been extensively studied in numerical experiments and has also been the subject of some articles concerning its analytical properties (see  Fabrie and Galusinski \cite{MR1868354}, Plech{\'a}{\v{c}} and {\v{S}}ver{\'a}k \cite{MR2012858}, Temam  \cite{MR0237972}). However, only recently Rusin \cite{rusin} proved existence of global weak solutions in $\R^3$.

\bt[Thm. 4.2,  Rusin \cite{rusin}] For any $\eps > 0$ and $\ueps _0 \in L^2 (\R^3)$, there exists a weak solution to (\ref{eqn:rusin}).
\et

\br \, Rusin \cite{rusin} also proved that when $\eps$ goes to zero, solutions to (\ref{eqn:rusin})  converge in $L^3 _{loc} (\R^3 \times \R_+) $ to a suitable (in the sense of Caffarelli-Kohn-Nirenberg \cite{MR673830}) solution to the Navier-Stokes equation.
\er

\medskip

\noindent Note that  since 
\be
\label{eqn:u-times-unablau}
\int _{\R^3} \ueps \left( \ueps \cdot \nabla \right) \ueps \, dx = - \frac{1}{2} \int _{\R^3} |\ueps|^2 div \, \ueps \, dx,
\ee
the nonlinear part vanishes when we multiply (\ref{eqn:rusin}) by $\ueps$ and integrate. Then, as the linear part of (\ref{eqn:rusin}) fits in the framework of  Section \ref{section-linear-part} (see Example \ref{example-rusin}), we obtain an inequality of the form (\ref{eqn:ineq-fourier-splitting}), with $\a = 1$ and we can use the Fourier Splitting method and the results in Sections \ref{section-decay-l2-linear-part} and \ref{section-decay-hs-solution} to study decay of solutions to (\ref{eqn:rusin}).

\subsubsection{$L^2$ Decay} We first prove the following result:

\bt
\label{decay-rusin-l2}
Let $\eps > 0$, $\ueps _0 \in L^2 (\R^3)$ with $r ^{\ast} = r ^{\ast} (\ueps _0)$, such that $- \frac{3}{2} < r^{\ast} < \infty$. Then for a weak solution $\ueps$ to (\ref{eqn:rusin}) we have that

\begin{displaymath}
\Vert \ueps (t) \Vert _{L^2} ^2 \leq C (1 + t) ^{- \min \{ \frac{3}{2} + r^{\ast}, \frac{5}{2}\}}.
\end{displaymath}

\et

{\bf Proof:}  As before, we proceed formally, by using the Fourier Splitting method for solutions which we assume are regular enough. The estimate for weak solutions is obtained as in Theorem \ref{decay-l2-qg}. Let  $B(t) = \{\xi \in \R^3: |\xi|^2 \leq g(t) \}$, for a nonincreasing, continuous $g$ such that $g(0) = 1$. By example (\ref{example-rusin}) and using the Fourier Splitting method as in Theorem \ref{decay-linear-part-qg} we obtain

\begin{align}
\label{eqn:main-estimate-fs}
\frac{d}{dt} \left( \exp \left( \int _0 ^t 2 g^2 (s) \, ds \right)  \Vert \ueps (t) \Vert _{L^2} ^2 \right) \leq \;\;\;\;\;\;\;\;\;\;\;\;\;\;\;\; \\ g^2 (t) \left( \exp \left( \int _0 ^t 2 g^2 (s) \, ds \right) \right) \int _{B(t)} |\widehat{\ueps} (\xi, t)| ^2 \, d \xi. \notag
\end{align}
We now have

\begin{displaymath}
\widehat{\ueps} (\xi, t) = e^{t \mathcal{M} (\xi)} \widehat{\ueps _0} (\xi)- \int _0 ^t e^{(t - s) \mathcal{M} (\xi)} G(\xi, s) \, ds
\end{displaymath}
where $e^{t \mathcal{M} (\xi)}$ is as in (\ref{eqn:sol-fund-frec}) and, denoting the Fourier transform by $\mathcal{F}$, 
\begin{displaymath}
G(\xi, s) = {\mathcal F} \left( (\ueps \cdot \nabla) \ueps + \frac{1}{2} (div \, \ueps)  \ueps \right).
\end{displaymath}
As 

\be
\label{eqn:equality-unablau}
(\ueps \cdot \nabla) \ueps = \nabla \cdot (\ueps \otimes \ueps) - (div \, \ueps)  \ueps
\ee
and 

\begin{displaymath}
\left| \widehat{(div \, \ueps)  \ueps} (\xi, t)\right| \leq \Vert div \, \ueps (t) \Vert _{L^2}  \Vert \ueps(t) \Vert _{L^2},
\end{displaymath}
we obtain

\begin{displaymath}
|G (\xi, t)| \leq |\xi| \Vert \ueps(t) \Vert _{L^2} ^2 + \frac{1}{2} \Vert div \, \ueps (t) \Vert _{L^2}  \Vert \ueps (t) \Vert _{L^2} \leq C |\xi| \Vert \ueps (t) \Vert _{L^2} ^2.
\end{displaymath}
Thus,
\begin{eqnarray*}
\left| \int _0 ^t e^{(t -s) \mathcal{M} (\xi)} G(\xi,s) \, ds \right| & \leq & C \int _0 ^t e^{- C (t - s) |\xi|^2} \, |\xi| \Vert \ueps(s) \Vert _{L^2} ^2 \, ds \nonumber \\ & \leq & C |\xi| \left( \int _0 ^t \Vert \ueps(s) \Vert _{L^2} ^2 \, ds \right).
\end{eqnarray*}
 Suppose now that $\Vert \ueps(t) \Vert _{L^2} ^2 \leq C (1 + t) ^{- \beta}$, for  some $0 \leq \beta$. We then have

\begin{displaymath}
\int _{B(t)} \left( \int _0 ^t e^{(t -s) \mathcal{M} (\xi)} G(\xi,s) \, ds \right) ^2 d \xi \leq C |\xi|^5 (1 + t) ^{2 (1 - \beta)},
\end{displaymath}
which leads, after choosing $g^2 (t) = \a (1 +t) ^{-1}$ and for large enough $\a > 0$, to

\begin{eqnarray}
\label{eqn:int-ball}
\int _{B(t)} |\widehat{\ueps} (\xi, t)|^2 \, d \xi & \leq & C \int _{B(t)} |e^{t \mathcal{M} (\xi)} \widehat{\ueps _0}|^2 \, d \xi + C \int _{B(t)} \left( \int_0 ^t e^{(t -s) \mathcal{M} (\xi)} G(\xi,s) \, ds \right) ^2 \, d \xi \nonumber \\ & \leq & C \Vert  e^{t \mathcal{M} (\xi)} \widehat{\ueps _0} \Vert _{L^2} ^2 + C g^5 (t) (1 + t) ^{2 (1 - \beta)} \nonumber \\ & \leq & C (t + 1) ^{- \left( \frac{3}{2} + r^{\ast} \right)} + C (1 + t) ^{- \left( \frac{1}{2} +2 \beta \right)} \nonumber \\ & \leq & C (t + 1) ^{- \min\{\frac{1}{2}+2\b,\frac{3}{2} + r^{\ast}\}},
\end{eqnarray}
where we used Theorem \ref{characterization-decay-l2} for the decay of the linear part. From (\ref{eqn:main-estimate-fs}), (\ref{eqn:int-ball}) and our choice of $g$ we obtain

\be
\label{eqn:estimate-cases}
\frac{d}{dt} \left( (t + 1) ^{\a}  \Vert \ueps (t) \Vert _{L^2} ^2 \right) \leq C (t + 1) ^{\a - 1} (t + 1) ^{- \min\{\frac{1}{2}+2\b,\frac{3}{2} + r^{\ast}\}}.
\ee
We start with $\b=0$, i. e. the known estimate  $\Vert \ueps(t) \Vert _{L^2} ^2 \leq C$. Then we split the study of $(\ref{eqn:estimate-cases})$ in the two cases  $\frac{3}{2} +r^{\ast} \leq   \frac{1}{2}$ and $  \frac{1}{2}\leq \frac{3}{2} +r^{\ast}$. In the first case, i.e. when $r^{\ast} \leq -1$, we obtain

\be
\label{eqn:decay-as-linear}
\Vert \ueps (t) \Vert _{L^2} ^2 \leq C  (t + 1) ^{- \left( \frac{3}{2} + r^{\ast} \right)}.
\ee
In the second case, we obtain 

\begin{displaymath}
\Vert \ueps (t) \Vert _{L^2} ^2 \leq C  (t + 1) ^{- \frac{1}{2}},
\end{displaymath}
i.e. we bootstrapped to $\beta = \frac{1}{2}$. Using this in $(\ref{eqn:estimate-cases})$, we see we have to separate again the study in two cases, $ \frac{3}{2} + r^{\ast} \leq \frac{3}{2}$ and $\frac{3}{2} \leq  \frac{3}{2} + r^{\ast}$. In the first case, i.e. when $r^{\ast} \leq 0$, once more we obtain (\ref{eqn:decay-as-linear}). In the second situation, i.e. when $r^{\ast} \geq 0$, we have improved to $\beta = \frac{3}{2}$. But then

\begin{displaymath}
\int _0 ^t \Vert \ueps (s) \Vert _{L^2} ^2 \, ds \leq C,
\end{displaymath}
so 

\begin{displaymath}
\int _{B(t)} \left( \int_0 ^t e^{(t -s) \mathcal{M} (\xi)} G(\xi,s) \, ds \right) ^2 \, d \xi \leq C g^5 (t).
\end{displaymath}
Then (\ref{eqn:int-ball}) becomes

\begin{eqnarray*}
\int _{B(t)} |\widehat{\ueps} (\xi, t)|^2 \, d \xi & \leq & C \int _{B(t)} |e^{t \mathcal{M} (\xi)} \widehat{\ueps _0}|^2 \, d \xi + C \int _{B(t)} \left( \int_0 ^t e^{(t -s) \mathcal{M} (\xi)} G(\xi,s) \, ds \right) ^2 \, d \xi \\ & \leq & C \Vert  e^{t \mathcal{M} (\xi)} \widehat{\ueps _0} \Vert _{L^2} ^2 + C g^5 (t)  \\ & \leq & C (t + 1) ^{- \left( \frac{3}{2} + r^{\ast} \right)} + C (1 + t) ^{- \frac{5}{2}} \leq C (t + 1) ^{- \min \{\frac{3}{2} + r^{\ast}, \frac{5}{2} \} }.
\end{eqnarray*}
Using this in (\ref{eqn:main-estimate-fs}) yields the conclusion of the proof. $\Box$

\smallskip 

We now address the decay of $w = \ueps - \bar{u}$, where $\bar{u} (x, t) = e^{t \mathcal{L}} u_0$ and  $\mathcal{L}$ is as defined in  (\ref{example-rusin}). 

\bt 
\label{decay-nonlinear-part-rusin}
Let $\eps > 0$, $\ueps _0 \in L^2 (\R^3)$, and $r^{\ast} = r^{\ast} (u_0)$ with $- \frac{3}{2} < r^{\ast} < \infty$. Then

\begin{displaymath}
\Vert \ueps (t) - \bar{u} (t) \Vert _{L^2} ^2 \leq C (1 + t) ^{- \min \{ \frac{7}{4} + r^{\ast}, \frac{5}{2} \}}
\end{displaymath}

\et

{\bf Proof:} The difference  $w$ solves the equation

\begin{displaymath}
w_t + (\ueps \cdot \nabla) \ueps + \frac{1}{2} (div \, \ueps)  \ueps = \Delta w + \nabla div \, w.
\end{displaymath}
Multiplying by $w$, integrating and using (\ref{eqn:u-times-unablau}) and (\ref{eqn:equality-unablau}) yields

\begin{eqnarray}
\label{eqn:energy-estimate-w}
\frac{d}{dt} \Vert w(t) \Vert _{L^2} ^2 & + & 2 \Vert \nabla w (t) \Vert _{L^2} ^2 + 2 \Vert div \, w (t) \Vert _{L^2} ^2 = \nonumber \\  & - & 2 \int _{\R^3} \bar{u} \left( (\ueps \cdot \nabla) \ueps + \frac{1}{2} (div \, \ueps) \ueps \right) \, dx.
\end{eqnarray}
Let $B(t) = \{ \xi \in \R^3: |\xi| \leq g(t) \}$, for a positive, decreasing to be determined function $g$ with $g(0) = 1$. As usual in the Fourier Splitting method

\begin{eqnarray}
\label{eqn:inequality-bcomp}
\Vert \nabla w (t) \Vert _{L^2} ^2 & + & \Vert div \, w \Vert _{L^2} ^2 \geq 2 g^2 (t) \int _{B(t) ^c} |\widehat{w} (\xi,t)| ^2 \, d \xi \nonumber \\ & = & 2 g^2 (t) \int _{\R^3} |\widehat{w} (\xi,t)| ^2 \, d \xi  - 2 g^2(t) \int _{B(t)} |\widehat{w} (\xi,t)| ^2 \, d \xi .
\end{eqnarray}
Now, as 

\begin{eqnarray*}
 \int_{\R^3} \bar{u} \, (div \, \ueps) \ueps \, dx  & = &  \int _{\R^3} \widehat{\bar{u}} \, \widehat{(div \, \ueps) \ueps} \, d \xi  \\ & \leq&  C \int_{\R^3} |\xi| |\widehat{\bar{u}}  \left( \widehat{\ueps} \ast \widehat{\ueps} \right)| \, d \xi  \leq  C \int _{\R^3} | \nabla \bar{u} (\ueps \otimes \ueps) | \, dx 
\end{eqnarray*}
we have that

\begin{eqnarray}
\label{eqn:inequality-nonlinear-term}
\left| \int _{\R^3} \bar{u} \, \left( (\ueps \cdot \nabla) \ueps + \frac{1}{2} (div \, \ueps) \ueps \right) \, dx \right| & = & \left| \int _{\R^3} \bar{u} \, \left( \nabla \cdot (\ueps \otimes \ueps) - \frac{1}{2}  (div \,\ueps)  \ueps \right) \, dx \right| \nonumber \\ & \leq &   C  \int _{\R^3}\left| \nabla \bar{u} \, (\ueps \otimes \ueps)\right|  \, dx \nonumber \\ & \leq & C \Vert \nabla \bar{u} (t)  \Vert _{L^{\infty}} \Vert \ueps (t) \Vert _{L^2} ^2  \nonumber \\ & \leq & C (1 + t)^{- \frac{5}{4}} \Vert \ueps (t) \Vert _{L^2} ^2 
\end{eqnarray}
where we used for $\nabla \bar{u}$ an analog to a standard estimate for solutions to the heat equation (see (2.3'), page 474 in Kato \cite{Kato1984}). Using (\ref{eqn:inequality-bcomp}) and (\ref{eqn:inequality-nonlinear-term}) in (\ref{eqn:energy-estimate-w}) we obtain, after multiplying by the appropriate integrating factor

\begin{align}
\label{eqn:big-ineaquality}
\frac{d}{dt} \left( exp \left( \int _0 ^t g^2(s) \, ds \right) \Vert w(t) \Vert _{L^2} ^2 \right)  \leq  \;\;\;\;\;\;\;\;\;\;\;\;\;\;\;\; \\  exp \left( \int _0 ^t g^2(s) \, ds \right)\left( g^2(t)\int _{B(t)} |\widehat{w} (\xi,t)| ^2 \, d \xi  +  C (1 + t)^{- \frac{5}{4}} \Vert \ueps (t) \Vert _{L^2} ^2  \right).\notag
\end{align}
As $w_0 = 0$, 

\begin{displaymath}
\widehat{w} (\xi,t) = \int_0 ^t e^{(t - s) \mathcal{M} (\xi)} \mathcal{F} \left( \nabla \cdot (u \otimes u) - \frac{1}{2}(div \,u)  u \right) (\xi, s)  \, ds
\end{displaymath}
leads to

\begin{displaymath}
|\widehat{w} (\xi, t)| \leq C \int_0 ^t e^{(t - s) \mathcal{M} (\xi)} |\xi| |\widehat{\ueps} \ast \widehat{\ueps} (\xi)| \, d \xi
\end{displaymath}
and hence
\be
\label{eqn:estimate-in-ball}
\int _{B(t)} |\widehat{w} (\xi,t)|^2 \, d \xi \leq g^5 (t) \left( \int _0 ^t \Vert \ueps(t) \Vert _{L^2} ^2 \, ds\right) ^2.
\ee
Let $g^2(t) = \a (1 + t) ^{-1}$, with $\a >0$ large enough, and using (\ref{eqn:estimate-in-ball}) in (\ref{eqn:big-ineaquality}) yields\begin{eqnarray}
\label{eqn:final-eqn-decay-nonlinear}
\frac{d}{dt} \left( (t + 1) ^{\a} \Vert w(t) \Vert _{L^2} ^2 \right) & \leq & C (t + 1) ^{\a - \frac{7}{2}} \left( \int _0 ^t \Vert \ueps(s) \Vert _{L^2} ^2 \, ds \right) ^2 \nonumber \\ & + &  C (1 + t)^{\a - \frac{5}{4}} \Vert \ueps (s) \Vert _{L^2} ^2.
\end{eqnarray}
Now we  use the decay rates  obtained in Theorem \ref{decay-rusin-l2}. Suppose first that $r^{\ast} = - \frac{1}{2}$. Then

\begin{displaymath}
\int _0 ^t \Vert \ueps (s) \Vert _{L^2} ^2 \, ds = \ln (1 + t),
\end{displaymath}
so (\ref{eqn:final-eqn-decay-nonlinear}) leads to

\begin{displaymath}
\frac{d}{dt} \left( (t + 1) ^{\a} \Vert w(t) \Vert _{L^2} ^2 \right)  \leq  C (t + 1) ^{\a - \frac{7}{2}} \ln ^2 (1 + t)  +   C (1 + t)^{\a - \frac{9}{4}}.
\end{displaymath}
Since $ \ln(t+1)^2(t+1)^{-\frac{5}{2}} \leq (t+1)^{-\frac{5}{4}}$, a simple computation leads to

\be \label{r1/2}
\|w(t)\|_2^2=\Vert \ueps (t) - \bar{u} (t) \Vert _{L^2} ^2 \leq C (1 + t) ^{- \frac{1}{4}}
\ee

Now assume $r^{\ast} \neq - \frac{1}{2}$. 

\noindent {\em Case 1.} If $\min \{ \frac{3}{2} + r^{\ast}, \frac{5}{2} \} = \frac{3}{2} + r^{\ast}$, then we have two possibilities, i.e. $\frac{3}{2} + r^{\ast} <1$ and $\frac{3}{2} + r^{\ast} >1$. In the first case, as 

\begin{displaymath}
\int _0 ^t \Vert \ueps(t) \Vert _{L^2} ^2 \, ds \leq C(t+1)^{-(\frac{1}{2}+r^*)}
\end{displaymath}
from (\ref{eqn:final-eqn-decay-nonlinear}) we obtain

\begin{displaymath}
\frac{d}{dt} \left( (t + 1) ^{\a} \Vert w(t) \Vert _{L^2} ^2 \right)  \leq  C (t + 1) ^{\a - \frac{9}{2} - 2 r^{\ast}}  +   C (1 + t)^{\a - \frac{11}{4} - r^{\ast}}.
\end{displaymath}
After integrating and comparing exponents we obtain
\be \label{2i}
\|w(t)\|_{L^2} ^2= \Vert \ueps (t) - \bar{u} (t) \Vert _{L^2} ^2 \leq C (1 + t) ^{- \left( \frac{7}{4} + r^{\ast} \right)}.
\ee
In the second case, $1< r^* +\frac{3}{2} \leq \frac{5}{2} $ or equivalently $-\frac{1}{2} \leq r^*\leq 1$. Then, from (\ref{eqn:final-eqn-decay-nonlinear}) we have

\be \label{help}
\frac{d}{dt} \left( (t + 1) ^{\a} \Vert w(t) \Vert _{L^2} ^2 \right)  \leq  C (t + 1) ^{\a - \frac{7}{2} } +   C (1 + t)^{\a - \frac{11}{4} - r^{\ast}}.
\ee
First, if in (\ref{help}) we use $ r^* +\frac{11}{4} \leq \frac{7}{2}$,  hence $-\frac{1}{2}< r^* \leq \frac{3}{4}$, after integrating we obtain 

\be \label{a}
\|w(t)\|_{L^2} ^2=\Vert \ueps (t) - \bar{u} (t) \Vert _{L^2} ^2 \leq C (1 + t) ^{- \frac{7}{4}-r^*}.
\ee
If $ \frac{7}{2}\leq r^* +\frac{11}{4} $, hence $\frac{3}{4}\leq r^* \leq 1$, the same method leads to 

\be \label{b}
\|w(t)\|_{L^2} ^2=\Vert \ueps (t) - \bar{u} (t) \Vert _{L^2} ^2 \leq C (1 + t) ^{- \frac{5}{2}}.
\ee
\noindent {\em Case 2.} If $\min \{ \frac{3}{2} + r^{\ast}, \frac{5}{2} \} = \frac{5}{2}$, then $r^* \geq 1$. In this case (\ref{eqn:final-eqn-decay-nonlinear}) yields

\begin{displaymath}
\frac{d}{dt} \left( (t + 1) ^{\a} \Vert w(t) \Vert _{L^2} ^2 \right)  \leq  C (t + 1) ^{\a - \frac{7}{2} }  +   C (1 + t)^{\a - \frac{15}{4} - r^{\ast}}.
\end{displaymath}
After integrating we obtain 

\be \label{2ii}
\|w(t)\|_{L^2} ^2=\Vert \ueps (t) - \bar{u} (t) \Vert _{L^2} ^2 \leq C (1 + t) ^{- \frac{5}{2}}.
\ee
Combining the estimates obtained in (\ref{r1/2}), (\ref{2i}), (\ref{a}),(\ref{b}) and (\ref{2ii}) we prove Theorem. $\Box$

As in the case of the dissipative quasi-geosotrophic equation, we obtain lower bounds for the decay of $\ueps$ along the lines of Theorem \ref{lower-bounds-qg}, using the decays from Theorems \ref{decay-rusin-l2} and \ref{decay-nonlinear-part-rusin}. We omit this proof, as it consists of simple computations.

\bt 
\label{lower-bound-rusin}
Let $\ueps _0 \in L^2 (\R^3), r^{\ast} = r^{\ast} (u_0)$. Then for $-\frac{3}{2} <  r^* \leq 1$ we have that

\begin{displaymath}
\Vert \ueps (t) \Vert _{L^2} ^2 \geq C (1 + t) ^{- \left( \frac{3}{2} + r^{\ast} \right)}.
\end{displaymath}
\et
Now, combining the estimates from Theorems \ref{decay-rusin-l2} and \ref{lower-bound-rusin} we obtain the following result.

\bt
\label{summary-thm-rusin}
Let $\ueps _0 \in L^2 (\R^3), r^{\ast} = r^{\ast} (u_0)$. Then for $-\frac{3}{2} <  r^* \leq 1$, there exist $C_1, C_2 > 0$ such that

\begin{displaymath}
C_1 (1 + t) ^{- \left( \frac{3}{2} + r^{\ast} \right)} \leq \Vert \ueps (t) \Vert _{L^2} ^2 \leq C_2 (1 + t) ^{- \left( \frac{3}{2} + r^{\ast} \right)}.
\end{displaymath}
If $r^{\ast} > 1$, then

\begin{displaymath}
\Vert \ueps (t) \Vert _{L^2} ^2 \leq C (1 + t) ^{- \frac{5}{2}}.
\end{displaymath}
\et

\br 
\label{remark-comparison} \, The estimates for the compressible approximation (\ref{eqn:rusin}) summarized in Theorem \ref{summary-thm-rusin} are the same as those obtained in Theorem 6.5 in Bjorland and M.E. Schonbek \cite{MR2493562} for the Navier-Stokes equations. Note that the range of values for which the lower bound is valid in part $1.$ of that Theorem in \cite{MR2493562} can be extended to the same range as in  Theorem \ref{summary-thm-rusin} by using the arguments in Theorem \ref{lower-bounds-qg}. Thus, the stabilizing nonlinear damping term $\frac{1}{2} (div \, \ueps) \ueps$ provides enough dissipation so as to  have an energy inequality in (\ref{eqn:rusin}), but does not alter the range of values of $r^{\ast}$ for which the linear part has slower decay (cf. Proof of Theorem \ref{lower-bounds-qg}).  \er

\subsubsection{$\dot{H}^s$ Decay} The goal is to prove the following Theorem.

\bt 
\label{decay-hs-rusin}
Let $\ueps _0 \in H^s (\R^3), s \geq 1$, $r^{\ast} = r^{\ast} (\ueps _0)$. Then

\begin{displaymath}
\Vert \ueps (t) \Vert^2  _{\dot{H}^s} \leq C (1 + t) ^{- \left( s + \min \{ \frac{5}{2}, r^{\ast} + \frac{3}{2}\}\right)}.
\end{displaymath}
\et

{\bf Proof:} we follow ideas along the lines of those in the proof of Theorem 2.4 from M.E. Schonbek and T. Schonbek \cite{MR2001105}. We first prove the preliminary decay

\be
\label{eqn:preliminary-decay-hs-rusin}
\Vert \ueps (t) \Vert^2  _{\dot{H}^s} \leq C (1 + t) ^{- \min \{ \frac{5}{2}, r^{\ast} + \frac{3}{2}\}}.
\ee
We apply $\Lambda^s$ to (\ref{eqn:rusin}), then multiply it by $\Lambda ^s u$ and integrate in space to obtain

\begin{displaymath}
\frac{d}{dt} \Vert \Lambda^s \ueps (t) \Vert ^2 _{L^2} + C \Vert \Lambda^{s+1} \ueps (t) \Vert _{L^2} ^2 = - C \int _{\R^3} \Lambda ^s \ueps  \Lambda^s \left( \nabla \cdot (\ueps \otimes \ueps) - \frac{1}{2}  (div \,\ueps)  \ueps \right) \, dx,
\end{displaymath}
For the first term on the right hand side we have

\begin{eqnarray*}
\left| \int _{\R^3 } \Lambda^s \ueps \, \Lambda^s \left(  \nabla \cdot (\ueps \otimes \ueps) \right)  \, dx \right|  & \leq & \Vert \Lambda ^{s +1} \ueps \Vert  _{L^2} \Vert \Lambda^s (\ueps \otimes \ueps) \Vert _{L^2}  \leq \eta \Vert \Lambda ^{s +1} \ueps \Vert  _{L^2} ^2 \\ & + & C(\eta) \Vert \Lambda^s (\ueps \otimes \ueps) \Vert _{L^2} ^2.
\end{eqnarray*}
By the Calculus Inequality, with $ \frac{1}{2} = \frac{1}{p} + \frac{1}{q}$, for $1 < p,q < \infty$, we obtain

\be
\label{eqn:first-ineq-hs}
\Vert \Lambda^s (\ueps \otimes \ueps) (t) \Vert _{L^2} ^2 \leq C \Vert \ueps (t) \Vert _{L^q} ^2 \Vert \Lambda^s \ueps (t) \Vert _{L^p} ^2.
\ee
Interpolating between the $L^2$ and $H^s$ norms, we obtain a Maximum Principle $\Vert u (t) \Vert _{L^q} \leq \Vert u_0 \Vert _{L^q}$, which together with the Hardy-Littlewood-Sobolev inequality lead to

\be
\label{eqn:second-ineq-hs}
\Vert \Lambda^s \ueps (t) \Vert _{L^p} ^2 \leq  C \Vert \ueps (t) \Vert _{L^q} ^2 \Vert \Lambda ^{s + \frac{3}{2} \left( 1 - \frac{2}{p} \right) } \ueps (t) \Vert _{L^2} ^2 \leq C \Vert \Lambda ^{s + \frac{3}{2} \left( 1 - \frac{2}{p} \right) } \ueps (t) \Vert _{L^2} ^2.
\ee
Then

\begin{displaymath}
\left| \int _{\R^3 } \Lambda^s \ueps \, \Lambda^s \left(  \nabla \cdot (\ueps \otimes \ueps) \right)  \, dx \right| \leq \eta \Vert \Lambda ^{s +1} \ueps (t) \Vert  _{L^2} ^2 + C(\eta) \Vert \Lambda ^{s + \frac{3}{2} \left( 1 - \frac{2}{p} \right) } \ueps (t) \Vert _{L^2} ^2.
\end{displaymath}
Now for the second term

\begin{eqnarray*}
\left| \int _{\R^3 } \Lambda^s \ueps \, \Lambda^s (div \, \ueps)  \ueps \, dx \right| & \leq & \Vert \Lambda^{s + 1} \ueps (t) \Vert _{L^2} \Vert \Lambda ^{s - 1} (div \, \ueps (t) ) \Vert _{L^2} \\ & \leq &  \eta  \Vert \Lambda^{s + 1} \ueps (t)  \Vert _{L^2} ^2 + C(\eta) \Vert \Lambda ^{s - 1} (div \, \ueps (t)) \Vert _{L^2} ^2.
\end{eqnarray*}
and by y the Calculus Inequality again

\begin{displaymath}
\Vert \Lambda ^{s - 1} (div \, \ueps (t)) \Vert _{L^2} ^2 \leq \left( \Vert div \, \ueps (t) \Vert _{L^q} ^2 \Vert \Lambda ^{s - 1} \ueps (t) \Vert _{L^p} ^2 + \Vert \ueps (t) \Vert _{L^q} ^2 \Vert \Lambda ^{s - 1} div \, \ueps (t)  \Vert _{L^p} ^2 \right).
\end{displaymath}
We will now use the following Fractional Gagliardo-Nirenberg inequality.

\bt[Corollary 1.5, \cite{MR2883850}, see also {\tt arXiv.org/abs/1004.4287}.]
\label{gagliardo-nirenberg-fracionario}
Let $1 < p, p_0, p_1 < \infty$, $s, s_1 \in \R$, $0 \leq e \leq 1$. Then the fractional Gagliardo-Nirenberg inequality

\begin{displaymath}
\Vert u \Vert _{\dot{H}^s _p} \leq C \Vert u \Vert _{L^{p_0}} ^{1 - e} \Vert u \Vert _{\dot{H}^{s_1} _{p_1}} ^e
\end{displaymath}
is true in $\R^n$ if and only if
\begin{displaymath}
\frac{n}{p} - s = (1 - e) \frac{n}{p_0} + e \left( \frac{n}{p_1} - s_1 \right), \quad s \leq e \, s_1.
\end{displaymath}
\et
We then have

\begin{displaymath}
\Vert div \, \ueps (t) \Vert _{L^q} ^2 \leq C \Vert \ueps (t) \Vert _{L^2} ^{2 \left( 1 - \frac{1}{s +1} \left( 2 - \frac{2}{q} \right) \right) } \Vert \Lambda ^{s +1} \ueps (t) \Vert _{L^2} ^{\frac{1}{s +1} \left( 2 - \frac{2}{q} \right)} 
\end{displaymath}
and

\begin{displaymath}
\Vert \Lambda ^{s - 1} \ueps (t) \Vert _{L^p} ^2 \leq C \Vert \ueps (t) \Vert _{L^2} ^{2 \left( 1 - \frac{1}{s +1} \left( 2 - \frac{2}{p} \right) \right) } \Vert \Lambda ^{s +1} \ueps (t) \Vert _{L^2} ^{\frac{1}{s +1} \left( 2 - \frac{2}{p} \right)},
\end{displaymath}
which lead to 

\begin{eqnarray*}
\Vert div \, \ueps (t) \Vert _{L^q} ^2 \Vert \Lambda ^{s - 1} \ueps (t) \Vert _{L^p} ^2 & \leq &  C \Vert \ueps (t) \Vert _{L^2} \Vert \Lambda ^{s +1} \ueps (t) \Vert _{L^2} \\ & \leq &  \eta \Vert \ueps (t)  \Vert _{L^2} ^2 + C(\eta) \Vert \Lambda ^{s +1} \ueps (t) \Vert _{L^2} ^2.
\end{eqnarray*}
And using the Maximum Principle and the Hardy-Littlewood-Sobolev inequality again we obtain

\begin{eqnarray}
\label{eqn:ineq-bootstrap}
\Vert u \Vert _{L^q} ^2 \Vert \Lambda ^{s - 1} div \, \ueps (t)  \Vert _{L^p} ^2 & \leq &  C \Vert \ueps (t)  \Vert _{L^q} ^2 \Vert \Lambda ^{(s - 1) + \frac{3}{2} \left( 1 - \frac{2}{p} \right)} div \, \ueps (t) \Vert _{L^2} ^2 \nonumber \\ & \leq &  C \Vert \Lambda ^{s + \frac{3}{2} \left( 1 - \frac{2}{p} \right)} div \, \ueps (t)  \Vert _{L^2} ^2.
\end{eqnarray}
Now as

\begin{eqnarray*}
\frac{d}{dt} \Vert \Lambda^s \ueps (t)  \Vert _{L^2} ^2 & + & C \Vert \Lambda ^{s + 1} \ueps (t) \Vert _{L^2} ^2 \\ & \leq & \left| \int _{\R^3} \Lambda^s \ueps (t) \, \Lambda^s  \left( \nabla \cdot \left( \ueps (t) \otimes \ueps (t) \right)  -  \frac{1}{2}  (div \, \ueps(t))  \ueps (t) \right)  \, dx \right| \nonumber \\ & \leq & C_1(\eta) \Vert \Lambda ^{s +1} \ueps (t) \Vert _{L^2} ^2 + C_2(\eta) \Vert \Lambda ^{s + \frac{3}{2} \left( 1 - \frac{2}{p} \right)} \ueps (t) \Vert _{L^2} ^2 \nonumber \\ & + & C_3(\eta) \Vert \ueps (t) \Vert _{L^2} ^2
\end{eqnarray*} 
choosing a small enough $C_1(\eta)$ we obtain

\be
\label{eqn:third-ineq-hs}
\frac{d}{dt} \Vert \Lambda^s \ueps (t)  \Vert _{L^2} ^2 + C_1 \Vert \Lambda ^{s + 1} \ueps (t) \Vert _{L^2} ^2 \leq C_2 \Vert \Lambda ^{s + \frac{3}{2} \left( 1 - \frac{2}{p} \right)} \ueps(t)  \Vert _{L^2} ^2 + C_3 \Vert \ueps (t) \Vert _{L^2} ^2.
\ee
For $p < 6$, so that $\frac{3}{2} \left( 1 - \frac{2}{p} \right) < 1$, we have 

\begin{eqnarray*}
\Vert \Lambda ^{s + \frac{3}{2} \left( 1 - \frac{2}{p} \right)} \ueps (t)  \Vert _{L^2} ^2 & = & \int _{\R^3} |\xi| ^{2 \left( s + \frac{3}{2} \left( 1 - \frac{2}{p} \right) \right) } |\widehat{u}|^2 \, d \xi  \nonumber \\ & \leq & M^{2 \left( s + \frac{3}{2} \left( 1 - \frac{2}{p} \right) \right) } \int _{B_M} |\widehat{\ueps} (\xi,t) |^2 \, d \xi \nonumber \\ & +&  M^{- 2 \left( \frac{3}{p} - \frac{1}{2} \right)} \int _{B_M ^c} |\xi| ^{2 (s + 1)} |\widehat{\ueps} (\xi,t)|^2 \, d \xi \nonumber \\ & \leq &  M^{2 \left( s + \frac{3}{2} \left( 1 - \frac{2}{p} \right) \right) } \Vert \ueps (t) \Vert _{L^2} ^2 + M^{- 2 \left( \frac{3}{p} - \frac{1}{2} \right)} \Vert \Lambda ^{s + 1} \ueps (t) \Vert _{L^2} ^2
\end{eqnarray*} 
where $B_M = \{ \xi \in \R^3: |\xi|^2 \leq M \}$. Now for fixed $M$, we have that

\begin{displaymath}
\Vert \Lambda ^{s +1} \ueps (t) \Vert _{L^2} ^2 \geq \int _{B_M} |\xi| ^{2 (s +1)} |\widehat{\ueps} (\xi,t)|^2 \, d \xi \geq M \Vert \Lambda ^s \ueps (t) \Vert _{L^2} ^2 - M^{s +1} \Vert \ueps (t) \Vert _{L^2} ^2,
\end{displaymath}
so for large enough $M$

\begin{displaymath}
\frac{1}{2} \frac{d}{dt} \Vert \Lambda ^s \ueps (t) \Vert _{L^2} ^2 + C \Vert \Lambda ^{s + 1} \ueps (t) \Vert _{L^2} ^2 \leq C \Vert \ueps (t) \Vert _{L^2} ^2.
\end{displaymath}
Using the integrating factor $h(t) = e^{ct}$ and the decay from Theorem \ref{decay-rusin-l2} we prove the preliminary decay (\ref{eqn:preliminary-decay-hs-rusin}). To prove Theorem \ref{decay-hs-rusin}, we proceed as before. We first note that by interpolation we have that

\begin{displaymath}
\Vert \ueps (t) \Vert _{L^q} ^2 \leq C (1 + t) ^{- \min \{ \frac{5}{2}, r^{\ast} + \frac{3}{2} \}}.
\end{displaymath}
We then go back to (\ref{eqn:first-ineq-hs}) and (\ref{eqn:second-ineq-hs}), use this decay in (\ref{eqn:ineq-bootstrap}) and rewrite (\ref{eqn:third-ineq-hs}) as 

\begin{eqnarray*}
\frac{1}{2} \frac{d}{dt} \Vert \Lambda^s \ueps (t) \Vert _{L^2} ^2 + C \Vert \Lambda ^{s + 1} \ueps (t) \Vert _{L^2} ^2 & \leq & (1 + t)^{- \min \{ \frac{5}{2}, r^{\ast} + \frac{3}{2} \} } \Vert \Lambda ^{s + \frac{3}{2} \left( 1 - \frac{2}{p} \right)} \ueps (t)  \Vert _{L^2} ^2 \\ & + & (1 + t)^{- \min \{ \frac{5}{2}, r^{\ast} + \frac{3}{2} \} }.
\end{eqnarray*}
We finish the proof as in the case of the dissipative quasigeostrophic equation, see page 368 in M.E. Schonbek and T. Schonbek  \cite{MR2001105}. $\Box$

\bibliography{Cesar-MariaElena-March2014}

\bibliographystyle{plain}

\end{document}